\documentstyle{amsppt}
%
%
%
%

\def\N{\Bbb N }
\def\P{\Bbb P }

\def\Z{\Bbb Z }
\def\C{\Bbb C }
\loadeurm
\define\k{{\eurm k}}
\define\ce{{\eurm c}}
\define\ka{{\eurm k}}

\define\nl{\hfil\newline}
\define\ldot{\,.\,}
\redefine\l{\lambda}

\redefine\d{\delta}
\define\w{\omega}

\redefine\i{{\,\text{{\rm i}}\,}}

\define\mapleft#1{\smash{\mathop{\longleftarrow}\limits^{#1}}}
\define\mapright#1{\smash{\mathop{\longrightarrow}\limits^{#1}}}
\define\mapdown#1{\Big\downarrow\rlap{
   $\vcenter{\hbox{$\scriptstyle#1$}}$}}
\define\mapup#1{\Big\uparrow\rlap{
   $\vcenter{\hbox{$\scriptstyle#1$}}$}}
\define\cint #1{\frac 1{2\pi\i}\int_{C_{#1}}}
\define\cintt{\frac 1{2\pi\i}\int_{C_{\tau}}}

\define\cinttt{\frac 1{24\pi\i}\int_{C_{\tau}}}

\define\im{\text{Im\kern1.0pt }}
\define\re{\text{Re\kern1.0pt }}
\define\res{\text{res}}
\redefine\deg{\operatornamewithlimits{deg}}
\define\ord{\operatorname{ord}}

\define\fpz{\frac {d }{dz}}
\define\dzl{\,{dz}^\l}
\define\pfz#1{\frac {d#1}{dz}}
\define\K{\Cal K}
\define\U{\Cal U}
\redefine\O{\Cal O}
\define\He{\text{\rm H}^1}
\define\Ho{\text{\rm H}^0}
\define\A{\Cal A}
\redefine\L{\Cal L}
\redefine\D{\Cal D^{1}}
\define\KN {Krichever-Novikov}
\define\Pif {{P_{\infty}}}
\define\Uif {{U_{\infty}}}
\define\Uifs {{U_{\infty}^*}}
\define\KM {Kac-Moody}
\define\Fl{\Cal F^\lambda}
\define\Fln{\Cal F^\lambda_n}
\define\Fn#1{\Cal F^{#1}}
\define\gb{\overline{\frak g}}
\define\G{\overline{\frak g}}
\define\Gb{\overline{\frak g}}
\redefine\g{\frak g}
\define\Gh{\widehat{\frak g}}
\define\gh{\widehat{\frak g}}
\define\Ah{\widehat{\Cal A}}
\define\Lh{\widehat{\Cal L}}

\define\iN{i=1,\ldots,N}
\define\iNi{i=1,\ldots,N,\infty}
\define\pN{p=1,\ldots,N}

\define\de{\delta}
\define\kndual#1#2{\langle #1,#2\rangle}
\define \nord #1{:\mkern-5mu{#1}\mkern-5mu:}
\define \MgN{{\Cal M}_{g,N}} 
\define \MgNe{{\Cal M}_{g,N+1}} 
\define \mpt{(M,P_1,P_2,\ldots, P_N,\Pif)} 
\define \mpp{(M,P_1,P_2,\ldots, P_N)} 
\define \MgNn{{\Cal M}_{g,N}^{(1)}} 
\define \MgNen{{\Cal M}_{g,N+1}^{(1)}} 
\define \Mgo{{\Cal M}_{g,0}} 
\define \mptn{(M,P_1,P_2,\ldots, P_N,\Pif,z_1,\ldots,z_N,z_\infty)}
\define \mppn{(M,P_1,P_2,\ldots, P_N,z_1,\ldots,z_N)} 
\define \sinf{{\widehat{\sigma}}_\infty}
\define\Wt{\widetilde{W}}
\define\St{\widetilde{S}}
\define\Wn{W^{(1)}}
\define\Wtn{\widetilde{W}^{(1)}}
\define\btn{\tilde b^{(1)}}
\define\bt{\tilde b}
\define\bn{b^{(1)}}
%

\define\tOmega{\Tilde\Omega}
\define\tw{\Tilde\omega}
\define\hw{\hat\omega}
\define\s{\sigma}
\define \car{{\frak h}}    
\define \bor{{\frak b}}    
\define \nil{{\frak n}}    
\define \vp{{\varphi}}
\define\bh{\widehat{\frak b}}  
\define\bb{\overline{\frak b}}  
\define\Vh{\widehat V}
\define\KZ{Knizhnik-Zamolodchikov}
\define\ai{{\alpha(i)}}
\define\ak{{\alpha(k)}}
\define\aj{{\alpha(j)}}

\global\newcount\numsec\global\newcount\numfor
\global\newcount\numfig
\gdef\profonditastruttura{\dp\strutbox}
\def\senondefinito#1{\expandafter\ifx\csname#1\endcsname\relax}
\def\SIA #1,#2,#3 {\senondefinito{#1#2}
\expandafter\xdef\csname #1#2\endcsname{#3}\else
\write16{???? ma #1,#2 e' gia' stato definito !!!!} \fi}

\def\etichetta(#1){(\veroparagrafo.\veraformula)
\SIA e,#1,(\veroparagrafo.\veraformula)
 \global\advance\numfor by 1
 \write16{ EQ \equ(#1) == #1  }}

\def\letichetta(#1){\veroparagrafo.\veraformula
\SIA e,#1,{\veroparagrafo.\veraformula}
\global\advance\numfor by 1
 \write16{ Sta \equ(#1) == #1}}

\def\tetichetta(#1){\veroparagrafo.\veraformula 
\SIA e,#1,{(\veroparagrafo.\veraformula)}
\global\advance\numfor by 1
 \write16{ tag \equ(#1) == #1}}

\def \FU(#1)#2{\SIA fu,#1,#2 }

\def\etichettaa(#1){(A\veroparagrafo.\veraformula)
 \SIA e,#1,(A.\veroparagrafo.\veraformula)
\global\advance\numfor by 1
 \write16{ EQ \equ(#1) == #1  }}

\def\getichetta(#1){Fig. \verafigura
 \SIA e,#1,{\verafigura}
 \global\advance\numfig by 1
 \write16{ Fig. \equ(#1) ha simbolo  #1  }}

\newdimen\gwidth
\def\BOZZA{
\def\alato(##1){
 {\vtop to \profonditastruttura{\baselineskip
 \profonditastruttura\vss
 \rlap{\kern-\hsize\kern-1.2truecm{$\scriptstyle##1$}}}}}
\def\galato(##1){ \gwidth=\hsize \divide\gwidth by 2
 {\vtop to \profonditastruttura{\baselineskip
 \profonditastruttura\vss
 \rlap{\kern-\gwidth\kern-1.2truecm{$\scriptstyle##1$}}}}}
}
\def\alato(#1){}
\def\galato(#1){}
\def\veroparagrafo{\number\numsec}\def\veraformula{\number\numfor}
\def\verafigura{\number\numfig}
\def\Eq(#1){\eqno{\etichetta(#1)\alato(#1)}}
\def\teq(#1){\tag{\tetichetta(#1)\hskip-1.6truemm\alato(#1)}}
\def\eq(#1){\etichetta(#1)\alato(#1)}
\def\Eqa(#1){\eqno{\etichettaa(#1)\alato(#1)}}
\def\eqa(#1){\etichettaa(#1)\alato(#1)}
\def\eqv(#1){\senondefinito{fu#1}$\clubsuit$#1\else\csname fu#1\endcsname\fi}
\def\equ(#1){\senondefinito{e#1}\eqv(#1)\else\csname e#1\endcsname\fi}

\newcount\refCount
\def\newref#1 {\advance\refCount by 1
\expandafter\edef\csname#1\endcsname{\the\refCount}}

\newref rBPZ
\newref rBeTo  
\newref rBeRs  
\newref rBono  
\newref rBrem 
\newref rEgOo  
\newref rFeRs   
\newref rFeWi   
\newref rGrOr   
\newref rHAG   
\newref rHit   
\newref rIv    
\newref rKaRa 
\newref rKnZ  
\newref rKod  
\newref	rKNFa 
\newref	rKNFb 
\newref	rKNFc 
\newref rNov  
\newref rRDS  
\newref rSad
\newref rSRS 
\newref rSLa 
\newref rSLb 
\newref rSLc 
\newref rSDiss 
\newref rSDeg 
\newref rSCt 
\newref rSHab 
\newref rSSS 
\newref rSSpt 
\newref rSeg  
\newref rShea 
\newref rSheb 
\newref rSha 
\newref rShhw 
\newref rShns 
\newref rTUY 
\newref rUcft 
%


\def\KnZ{\number\rKnZ}  

%
%
%
\magnification=1200
\vsize=21.5truecm
\hsize=16truecm
\hoffset=0.5cm\voffset=0.2cm
\baselineskip=15pt plus 0.2pt
\parskip=8pt
\newif\iftes

\def\testno{\tesfalse}
%
\testno
%
\iftes
\def\stamp{\quad\the\day.\the\month.\the\year.}
\def\ter#1{\hskip -1.5cm\hbox to 1.5cm{{\bf #1\hfill}}}
\BOZZA
\else
\def\stamp{}
\def\ter#1{}
\fi
\NoBlackBoxes
\TagsOnRight
\hfill Mannheimer Manuskripte 236

\hfill math.QA/9812083
\vskip 2cm
\topmatter
\title
           The Wess-Zumino-Witten-Novikov theory,
              Knizhnik-Zamolodchikov equations,
               and Krichever-Novikov algebras, I
\endtitle
\rightheadtext{WZNW theory\qquad\stamp }
\leftheadtext{M. Schlichenmaier, O.K. Sheinman \qquad\stamp }
\author Martin Schlichenmaier
\footnote"*\quad"{Partially~supported~by~the~Volkswagen-Stiftung
(RiP-program~at~Oberwolfach)
and the DFG - RFBR Project
\line{Nr.~436RUS~113/276/0(R)(DFG),~96-01-00055G~(RFBR)
\hfill}}
and Oleg K. Sheinman\footnotemark"*\quad"
\endauthor
\address
Martin Schlichenmaier,
Department of Mathematics and Computer Science,
University of Mannheim
D-68131 Mannheim, Germany
\endaddress
\email
schlichenmaier\@math.uni-mannheim.de
\endemail
\address
Oleg K. Sheinman,
Independent University of Moscow,
Ul. Marshala Biryuzova 4, kor.1, kv.65
Moscow, Russia, 123298
\endaddress
\email
sheinman\@landau.ac.ru
\endemail
\date \stamp
\enddate
\dedicatory Dedicated to Prof. S.P.Novikov in honour of his 60-th
birthday
\enddedicatory
\abstract
Elements of  a global operator approach to the
WZWN theory for  compact Riemann surfaces of arbitrary genus $g$ are given.
Sheaves  of representations  of affine Krichever-Novikov
algebras over a dense open subset of the moduli space of Riemann surfaces
(respectively of smooth, projective complex curves)
with $N$ marked points are introduced.
It is shown that the tangent space of the moduli space
at an arbitrary moduli point is isomorphic to a certain subspace of the
Krichever-Novikov vector field algebra  given by the
data of the moduli point.
This subspace is complementary to the direct sum of the two
subspaces containing the vector fields which  vanish at
the marked points, respectively which are regular at
 a fixed reference point.
For each representation of the affine algebra
$3g-3+N$ equations $(\partial_k+T[e_k])\Phi=0$
are given, where the elements $\{e_k\}$ are a basis of the
subspace,
and $T$ is the
Sugawara representation of the centrally extended
vector field algebra. For genus zero one obtains
the Knizhnik-Zamolodchikov equations in this way. The coefficients of the
equations for genus one are found in terms of Weierstra\ss-$\sigma$
function.
\endabstract
\endtopmatter
%
%
%
%
%
%
%
\define\kzint{1}  
\define\kzkn{2}   
\define\kzrep{3}  
\define\kzgen{4}  
\define\kzrat{5}  
\define\kzell{6}  
\define\kzapp{7}  
%

%
%
%
%
\document
\newpage
\vskip 1cm
\head
            \kzint . Introduction
\endhead
\vskip 0.2cm
\numsec=\kzint   
\numfor=1          
In the development of two-dimensional
conformal field theory it was realized
that the following problem is of special importance:
construct a bundle
over the {\it moduli space of punctured Riemann surfaces}
\footnote{See Section \kzgen\ for the precise definition what
me mean by ``moduli space of punctured Riemann surfaces''.}
with actions of the
gauge and the conformal algebras,
and  with a projectively flat connection.
 This problem
originates in the well-known paper of Knizhnik and Zamolodchikov
\cite\rKnZ , where the case of genus zero is considered. To each
puncture they assigned an irreducible module of the affine Kac-Moody
algebra associated to the gauge algebra.
 Over each punctured Riemann sphere they took the tensor
product of these modules as a fibre of the bundle.
 Expectations of sections of the so
obtained bundle were interpreted by them as $N$-point correlation
functions, where $N$ is the number of punctures.
  Making use of the
Sugawara construction they obtained an action of the Virasoro
algebra in each fibre and then defined the connection as
$\frac \partial {\partial z_p} + L^{(p)}_{-1}$, where $p=1,\ldots,N$
is the index of the punctures, $z_p$ is a local coordinate in a
neighbourhood of the $p$-th puncture and $L^{(p)}_{-1}$ is the
corresponding (rescaled) Sugawara operator of degree $-1$.
 One of the most famous
results of the paper \cite\rKnZ\ is the derivation of
equations which have as solutions  the sections
 of degree zero which are horizontal with respect to the
connection. The equations are nowadays known as the
Knizhnik-Zamolodchikov (KZ) equations, the connection as
the KZ connection.
The KZ equations for $N$-point correlation functions look as follows:
  $$\bigg(k\frac{\partial}{\partial z_p}-
     \sum\limits_{r\ne p}\frac{t^a_pt^a_r}{z_p-z_r}\bigg)
    <\phi_1(z_1)\ldots\phi_N(z_N)>=0\ ,
                                    \qquad p=1,\ldots,N\ .
                                      \Eq(KZeqi)
  $$
Here $k$ denotes some constant depending on the Kac-Moody algebra and the
level of the representation.
The $t_i^a$ for $i=1,\ldots,N$ are the (anti-hermitian)
representation matrices
for the $a^{th}$ generator of the (finite-dimensional) gauge algebra
in the representation associated to the point $z_i$
(operating on the $i^{th}$ field).
In \equ(KZeqi) a summation over the index $a$ is assumed.

D.Bernard was the first who considered the case of positive genus.
He realized that ``almost all the special features of a given WZW
model are encoded inside its zero modes'' \cite{\rBeTo, p.81}, but
"The lack of a precise definition of the action of the zero modes
$J^a_{0;\ j}$ on the correlation functions is very bad" \cite{\rBeRs, p.146}.
In other words, he immediately met a problem of defining
an action of the zero modes of the current operator  on  correlation
functions. The way-out he proposed was to consider correlation
functions (as well as corresponding primary fields) as functions
of not only the punctures but of $g$ additional parameters which are
elements of the finite-dimensional Lie group $G$ ($g$ denotes the
genus of the Riemann surface). The additional parameters are called
{\it twists}. From the recent point of view what he really
constructed is a connection over the {\it moduli space of
representations} of the fundamental group of the punctured Riemann
surface in the group $G$.  The  many interesting ideas contained in
the papers \cite{\rBeTo},\cite{\rBeRs} stimulated a number of
articles \cite{\rFeRs},\cite{\rFeWi},\cite{\rIv} in which  higher
genus \KZ\ equations are interpreted as equations of horizontality
with respect to a projectively flat connection over the moduli space
of representations. In particular, in \cite{\rFeRs},\cite{\rFeWi} the
mathematical foundation of the theory and some of its mathematical
results are formulated with more precision.

Hitchin \cite\rHit\ proposed his own approach to the problem. He
obtained a projectively flat connection on the moduli space of
representations by means of geometric quantization of certain integrable
systems.

One more direction was started by the important paper of Tsuchiya,
Ueno, Yamada \cite\rTUY\ . In contrary to the above mentioned works
(except the approach of \KZ) they constructed a projectively flat connection on
the space of punctured Riemann surfaces (more precisely, on  the
moduli space of stable curves with marked points).
In a most direct and consecutive form
their approach  used
the idea of deforming the complex structure on a Riemann surface by means
of cutting and twisting little circles around certain points.
They chose local coordinates (resp.  formal
neighbourhoods) at the punctures and then reproduced the \KZ\ construction
with respect to the chosen coordinates. To each puncture they assigned a
copy of the Kac--Moody loop algebra and a
copy of the Virasoro algebra connected with the
fixed local coordinate in the neighbourhood of that puncture. They
considered as  algebra of gauge symmetries of the theory
a central extension of the  direct sum of the local loop algebras
and as algebra of conformal symmetries the direct sum of the local
Virasoro algebras.
Starting from the local representations a
sheaf of representations was constructed over the moduli space.
The {\it conformal blocks} were introduced as a certain quotient sheaf.
It turned out to be a finite-dimensional vector bundle over the moduli
space. The above mentioned connection was constructed in this bundle.
Note that they were able to supply with their techniques a
mathematical proof of the Verlinde formula.
See also \cite{\rUcft} for a pedagogical presentation of this
approach.

It was clear from the very beginning that the basic objects of the
Wess-Zumino-Witten-Novikov (WZWN)
 theory on a Riemann surface are of  global nature. In fact, the
current and the energy-momentum tensor are abelian differentials (of
1-st and  2-d order respectively). The algebra annihilating
{conformal blocks} consists of (Lie algebra valued) meromorphic
functions with certain  polar behaviour \cite\rFeWi.  Nevertheless, most
authors used the above mentioned local approach
to gauge and conformal symmetries
due to Tsuchiya, Ueno and Yamada.
In \cite{\rKNFa-\rKNFc} Krichever and Novikov introduced their
generalizations of affine Kac-Moody and Virasoro algebras and pointed out a
new (global) treatment of gauge and conformal symmetries in
two-dimensional conformal field theory.  \KN\ algebras were studied
and generalized in
\cite{\rBrem},\cite{\rSLa-\rSHab},\cite{\rShea-\rShns}.
In \cite{\rBono},\cite{\rSSS}
the  Sugawara
construction for these algebras was considered.

The starting point of this article is our observation that it is
natural to consider multipoint \KN\ algebras as algebras of gauge
and conformal symmetries in WZWN theory for arbitrary genus
\footnote{Also I. Krichever expressed this idea to one of the
authors in a private talk.}.
These algebras have a genuine connection
with the basic geometrical objects of the theory, namely with
Riemann surfaces with punctures. There is no
problem of the action of zero modes of the current operator
in the \KN\ set-up. Moreover, \KN\ vector fields
 have a natural connection  with
deformations of punctured Riemann surfaces.
As it is shown below (Theorem~\kzgen.5) they can be related to
the Kuranishi tangent space of the moduli space of
punctured Riemann surfaces.
In a dual manner the deformations of punctured Riemann surfaces
can be described  by (meromorphic) quadratic
differentials with at most poles of order one at the punctures.
Note that the Krichever-Novikov duality (see Proposition~\kzkn.3)
supplies
such a dual description.
See also \cite{\rGrOr} for one more approach
for the two-point case.
\medskip

In Section \kzkn\  the necessary setup for the multipoint \KN\ algebras
in the situation of $N$ incoming points (corresponding to
the punctures $(P_1,P_2,\ldots,P_N)$  which can be moved)
and one outgoing point (corresponding to a fixed reference
point $\Pif$) are recalled.
The definition of  the \KN\ vector field algebras,
central extensions, affine algebras corresponding
to a fixed finite-dimensional Lie algebra, etc.
are given.  Necessary results are recalled.

In Section \kzrep\  Verma modules
for the higher genus multipoint affine algebras are introduced and studied.
The  Sugawara construction for these
\KN\ algebras is introduced along the lines of
\cite\rSSS .

In Section \kzgen\ the necessary moduli spaces of compact
Riemann surfaces with marked points are introduced.
As a technical tool we introduce also  enlarged moduli spaces,
with first order infinitesimal neighbourhoods at the punctures
as additional data.
Note that in this article we only deal
with the generic situation, i.e. the moduli point
corresponds to a generic curve and a generic
choice of marked points.
Hence, a priori our objects will only be defined over a dense
open subset of the moduli space.
Sheaf versions of the affine \KN\ algebras and the
Verma modules  are given.
The elements of the \KN\ vector field algebra define
tangent vectors along the moduli space.
In Theorem \kzgen.5 an explicit isomorphism of the
tangent space with a certain subspace of the vector field
algebra is given.
Let $\ \{X_k,\ k=1,\ldots,3g-3+N\}\ $ be a basis of the tangent
space and $l_k$
 the element of the vector field algebra
which  corresponds to $X_k$ under this isomorphism. For
sheaves of admissible representation
of the affine algebra (e.g. for the Verma modules)
the following set of
$3g-3+N$ equations
$$
\nabla_k\Phi:=\left(\partial_k+T[l_k]\right)\Phi=0\Eq(kzeint)
$$
is introduced as {\it formal KZ equations}.
Here $\partial_k$ is the derivation in direction of $X_k$ on the moduli
space and $T[l_k]$ the operator corresponding
under the Sugawara representation to the vector field
$l_k$.
Note that $T[l_k]$ operates vertically in the fibre.

In Section \kzrat\ the genus zero case is considered.
It is shown
how to obtain the original KZ equations
from \equ(kzeint).

In Section \kzell\ the genus one case is studied.
Explicit expressions for the coefficients in the
KZ equations in terms of Weierstra\ss-$\sigma$ function are derived.
In an appendix the KN basis for genus one and $N$ marked
points is given.

Altogether the proposed approach  enables us to
avoid some difficulties of the earlier approaches, to make use of
advantages of the global Sugawara construction, and to give a transparent
treatment of the geometric origin of
the coefficients in the Knizhnik-Zamolodchikov
equations.  It should be noted that even in the case $g=0$ our approach
provides a new derivation of the original KZ equations.

In the forthcoming part II of the article \cite{\rSSpt} we will study
the KZ equations on the sheaves of Verma modules in more detail,
develop further the structure theory and discuss conformal blocks and
projective flatness.


\vskip 1cm
\head
   \kzkn . The algebras of Krichever-Novikov type
\endhead
\numsec=\kzkn   
\def\kznum{\kzkn}
\numfor=1 
\vskip 0.2cm
\subhead
(a) The general set-up
\endsubhead

Let us recall here the set-up developed in \cite{\rSDiss},
\cite{\rSLa--\rSLc}.
Let $M$ be a compact Riemann surface of genus $g$, resp\. in the
language of algebraic geometry a smooth projective curve over $\C$.
Let
$$
I=(P_1,\ldots,P_N),\quad\text{and}\quad
O=(Q_1,\ldots,Q_L),\qquad  (N,L\ge 1)
$$
be disjoint tuples of ordered, distinct points (``marked points''
``punctures'') on the
curve. In particular, we assume $P_i\ne Q_j$ for every
pair $(i,j)$. The points in $I$ are
called the {\it in-points} the points in $O$ the {\it out-points}.
Let $A=I\cup O$ as a set.
In this article we are mainly dealing with $\# I=N\ge 1$ and
$\#O=1$.
Let $\rho$ be a  meromorphic differential on $M$, holomorphic on
$M^*:=M\setminus A$, with
positive residues at the points  in $I$, negative residues
at the points in $O$, and only purely imaginary periods.
By giving the residues (obeying the  condition
``sum over all residues equals zero'') there is a unique
such $\rho$.
If we choose an additional point $R\in M^*$ then the function
$$
u(P):=\re\int_R^P \rho
$$
is well-defined. Its level lines define a fibering
of $M^*$.
Every level line cuts the Riemann surface and separates the
in-points from the out-points.

Let $\K$ be the canonical line bundle.
Its associated sheaf of local sections is the sheaf of
holomorphic differentials.
Following the common practice we will usually not
distinguish between a line bundle and its associated invertible sheaf
(and even between the divisor class corresponding to a meromorphic
section of the line bundle).
For every $\l\in\Z$ we consider the bundle
$\ \K^\l:=\K^{\otimes \l}$. Here we use the usual convention:
$\K^0=\Cal O$ and $\K^{-1}=\K^*$
is the holomorphic tangent line bundle, (resp.
the sheaf  of holomorphic vector fields).
Indeed, after
fixing a theta characteristics, i.e. a bundle  $S$ with
$S^{\otimes 2}=\K$, it is possible to consider $\l\in \frac {1}{2}\Z$.
Denote by $\Fl$ the (infinite-dimensional) vector space of
global meromorphic sections  of $\K^\l$ which are holomorphic
on $M\setminus A$.
Special cases, which are of particular interest to us, are
the quadratic differentials ($\l=2$),
the  differentials ($\l=1$),
the functions  ($\l=0$), and
the vector fields ($\l=-1$).
The space of functions we will also denote by $\A$ and the
space of vector fields by $\L$.
By  multiplying  sections with functions
we again obtain sections. In this way
the space $\A$ becomes an associative algebra and the $\Fl$ become
$\A$-modules.

The vector fields in $\L$ operate on $\Fl$ by taking
the Lie derivative.
In local coordinates
 $$
{\nabla_e(g)}_|:=
L_e(g)_|:=(e(z)\fpz)\ldot (g(z)\dzl):=
\left( e(z)\pfz g(z)+\l\, g(z)\pfz e(z)\right)\dzl \ .\Eq(eLd)
$$
Here $e\in \L$ and $g\in \Fl$. To avoid cumbersome notation we
used the same symbol for the section and its representing
function.
If there is no danger of confusion we will do the same in the
following.
The space $\L$ becomes a Lie algebra with respect to \
the Lie bracket \equ(eLd)
and the $\Fl$ become Lie modules over $\L$.
Let us mention that by the action of $\L$ on $\A$ we
may define the
Lie algebra
$\D$ of differential operators of
degree $\le 1$ and the $\Fl$ become Lie modules over $\D$.
It is possible to extend these to differential operators of arbitrary
degrees, see \cite{\rSDiss,\rSCt,\rSHab} for
further information.
We will not need
this
additional structure here.
\definition{Definition \kznum.1}
The {\it Krichever-Novikov pairing} ({\it KN pairing}) is the
pairing between $\Fl$ and $\Fn {1-\l}$ given by
$$
\Fl\times\Fn {1-\l}\ \to\ \C,\qquad
\kndual {f}{g}:=\cintt f\cdot g\ ,\Eq(knpair)
$$
where $C_\tau$ is an arbitrary non-singular level line.
\enddefinition
Note that in \equ(knpair) the integral does not depend on
the level line chosen.
Using residues the pairing can  be described in a purely
algebraic manner as
$$
\kndual {f}{g}=\sum_{P\in I}\res_{P}(f\cdot g)=
-\sum_{Q\in O}\res_{Q}(f\cdot g)\ .\Eq(knpares)
$$


\subhead
(b) The almost-graded structure
\endsubhead

For the Riemann sphere ($g=0$) with quasi-global coordinate $z$
 and $I=(0)$ and $O=(\infty)$ the introduced vector field algebra is
the Witt algebra, i.e.  the algebra whose universal central extension
is the Virasoro algebra.
We denote for short this situation as the
{\it classical situation}.
Here it is of fundamental importance that
this algebra is a graded algebra.
For the higher genus case
(and for the multi-point situation for $g=0$) there is no such grading.
It was a fundamental observation by Krichever and Novikov
\cite{\rKNFa--\rKNFc} that a weaker concept, an
almost grading will do.
\definition{Definition \kznum.2}
(a) Let $\L$ be an (associative or Lie) algebra admitting a direct
decomposition as vector space $\ \L=\bigoplus_{n\in\Z} \L_n\ $.
$\L$ is called an {\it almost-graded}
({\it quasi-graded}, {\it generalized-graded})
algebra if (1) $\ \dim \L_n<\infty\ $ and (2)
there are constants $R$ and  $S$ with
$$
\L_n\cdot \L_m\quad \subseteq \bigoplus_{h=n+m-R}^{n+m+S} \L_h,
\qquad\forall n,m\in\Z\ .\Eq(eaga)
$$
The elements of $\L_n$ are called {\it homogeneous  elements of degree $n$}.
\nl
(b) Let $\L$ be an almost-graded  (associative or Lie) algebra
and $\Cal M$ an $\L$-module with
$\ \Cal M=\bigoplus_{n\in\Z} \Cal M_n\ $
as vector space. $\Cal M$ is called an {\it almost-graded}
({\it quasi-graded}, {\it generalized-graded}) module, if
(1) $\ \dim \Cal M_n<\infty\ $, and
(2) there are constants $R'$ and  $S'$ with
$$
 \L_m \cdot\Cal M_n\quad \subseteq \bigoplus_{h=n+m-R'}^{n+m+S'} \Cal M_h,
\qquad \forall n,m\in\Z\ .\qquad\Eq(egam)
$$
The elements of $\Cal M_n$ are called {\it homogeneous  elements of degree $n$}.
\enddefinition
\noindent
By a {\it weak almost grading} we understand an almost grading without
the requiring the finite-dimensionality of the homogeneous
subspaces.

For the 2-point situation, $I=\{P\}$ and $O=\{Q\}$, Krichever and Novikov
introduced an almost graded structure of the algebras and the modules
by exhibiting
special bases and defining their elements to be the
homogeneous elements.
By one of the authors its multi-point
generalization was given  \cite{\rSLc,\rSDiss}, again  by
exhibiting a special basis.
(See also Sadov \cite{\rSad} for some results in  similar directions.)
For every $n\in\Z$, $\pN$ a certain element
$f_{n,p}^\l\in\Fl$ is exhibited.
The $f_{n,p}^\l$ for $p=1,\ldots,N$ are a basis of a
subspace $\Fln$ and it is shown that
$$
\Fl=\bigoplus_{n\in\Z}\Fln\ .
$$
The subspace  $\Fln$ is called the {\it homogeneous subspace of degree $n$}.
\proclaim{Proposition  \kznum.3}\cite{\rSLc,\rSDiss}
(a) By the above definition the vector field algebra $\L$ and the
function algebra $\A$ are almost graded and the
modules $\Fl$ are  almost graded modules over them.
\nl
(b) The basis elements fulfil the duality relation
with respect to the KN pairing \equ(knpair)
$$\kndual {f_{n,p}^\l} {f_{m,r}^{1-\l}}=
 \cintt f_{n,p}^\l\cdot f_{m,r}^{1-\l}=\de_{-n}^{m}\cdot
\de_{p}^{r}\ , \Eq(edu)
$$
where $C_\tau $ is an arbitrary non-singular level line.
\endproclaim
\noindent
By \equ(edu) we see that the KN pairing is non-degenerate.
Let us introduce the following
notation:
$$
A_{n,p}:=f_{n,p}^0,\quad
e_{n,p}:=f_{n,p}^{-1},\quad
\w^{n,p}:=f_{-n,p}^1,\quad
\Omega^{n,p}:=f_{-n,p}^2 \ .\Eq(econc)
$$

The elements $f_{n,p}^\l$ have the following property
$$
\ord_{P_i}(f_{n,p}^\l)=(n+1-\l)-\d_i^p,\quad i=1,\ldots,N\ .
$$
The orders at the points in $O$ we will give
for the $(N,1)$ situation only (this is a short hand notation for
$\#I=N,\#O=1$).
Let us denote the single element in $O$ by $\Pif$

For $g=0$, or $g\ge 2$, $\l\ne 0,1$ and a generic
choice for the points in $A$
we have
$$
\ord_{\Pif}(f_{n,p}^\l)=-N\cdot(n+1-\l)
+(2\l-1)(g-1)\ .\Eq(ordi)
$$
By Riemann-Roch type arguments it is shown in \cite{\rSLa}
that there is up to a scalar multiple  only one such $f_{n,p}^\l$.
After choosing local coordinates $z_p$ at the points
$P_p$ the scalar may be fixed by requiring
$${f_{n,p}^\l}_|(z_p)=z_p^{n-\l}(1+O(z_p))\left(dz_p\right)^\l\ .$$
Due to the speciality of the occuring
divisors there is for a finite number of degrees $n$ a modified
prescription at the point $\Pif$ needed for the remaining cases.
In any case this is done without disturbing the orders at $I$ and
the KN duality.
The general description is given in \cite{\rSDiss, p.73}.
For the functions $A_{n,p}$ there are only modifications necessary if
all orders given by the generic rules are
nonpositive.
For the differentials $\w^{n,p}$ there are only modifications
necessary if all orders given by the generic rules are
nonnegative
or if there is exactly one pole of order 1 in the prescription.
In any case, modification  can only  appear for a finite number of
values of $n$.

To give an impression of the modification necessary we like to
give them for the $(1,1)$ situation
and $g\ge 2$.
For the 1-differentials modifications are
 necessary only for $-g\le n\le 0$.
The modified
basis elements of the differentials in
the $(1,1)$ situation and $g\ge 2$ are for $-g\le n\le -1$ $$
\ord_P(w^n)=-n-1,\quad
\ord_{\Pif}(w^n)=g+n,\quad
$$
and
$\
\ord_P(w^0)=-1,\quad
\ord_{\Pif}(w^0)=-1,\quad
\ $
with the additional condition that $\w^0$ has only imaginary periods.
In particular, $\ \w^0=\rho\ $ up to some multiplication with a scalar.
For the functions we set
for $-g\le n\le -1$
$$
\ord_P(A_n)=n,\quad
\ord_{\Pif}(A_n)=-g-n-1,\quad
$$
and $\ A_0=1\ $.
To fix the elements $A_n$ for $-g\le n\le -1$ we add suitable
multiples of $A_0$ such that the duality
$\ \langle w^0,A_n\rangle=0\ $
is fulfilled.
The $g=1$ situation will be covered in Section \kzapp.


For the basis elements $f_{n,p}^\l$ explicit
descriptions in terms of rational functions (for $g=0$),
the Weierstra\ss\ $\sigma$-function (for $g=1$), and prime forms
and theta functions (for $g\ge1$) are given in
\cite{\rSLb}.
For a description using Weierstra\ss\ $\wp$-function, see
\cite{\rRDS}, \cite{\rSDeg}.
The existence of such  an description is necessary in
our context because we want
to consider the above algebras  and modules over the
configuration space, resp\. the moduli space of pointed curves.
In particular, by the explicit representation one sees that
the elements vary ``analytically'' when
the complex structure of the Riemann surface is deformed

For further reference and as an illustration
 let us write down the basis elements
for $g=0$.
We choose a quasi-global coordinate $z$ such that
the point $\Pif$ is given by $z=\infty$.
Let the points $P_i$ be given by $z=z_i$ for $\iN$.
Clearly,
$$
f_{n,p}^\l(z)=
(z-z_p)^{n-\l}\bigg(\prod_{i=1\atop i\ne p}^N(z-z_i)\bigg)^{n-\l+1}
\bigg(\prod_{i=1\atop i\ne p}^N(z_p-z_i)\bigg)^{-n+\l-1}
dz^\l\ .\Eq(ebgo)
$$

An explicit description for $g=1$ can be found in Section \kzapp.

The constructed basis coincide with  the Virasoro basis
in the classical situation, and with  the
basis for the two-point situation in higher genus given
by Krichever and Novikov
\cite{\rKNFa--\rKNFc} (up to some index shift).

We need a finer description of the almost graded structure.
The following is shown in
\cite{\rSLc,\rSDiss}
\proclaim{Proposition  \kznum.4}
There exists constants $\ K,L\in\N\ $ such that
for all $n,m\in\Z$
$$
\align
A_{n,p}\cdot A_{m,r}&=\d_p^r\,A_{n+m,p}+
\sum_{h=n+m+1}^{n+m+K}\sum_{s=1}^N\alpha_{(n,p),(m,r)}^{(h,s)} A_{h,s}\ ,
\\
[e_{n,p},e_{m,r}]&=\d_p^r\,(m-n)\,e_{n+m,p}+
\sum_{h=n+m+1}^{n+m+L}\sum_{s=1}^N\gamma_{(n,p),(m,r)}^{(h,s)} e_{h,s}\ ,
\endalign
$$
with suitable coefficients $\ \alpha_{(n,p),(m,r)}^{(h,s)},
\gamma_{(n,p),(m,r)}^{(h,s)}\in\C$.
\endproclaim
The constants $K$ and $L$ can be explicitly calculated.
They depend on the genus $g$ and on the number of points $N$.
Again we give here only the result for the $(N,1)$ situation and $g\ne 1$.
$$
\aligned
L&=\cases 3g,&\text{$g\ne1$ and  $N=1$}
\\   3+\left[\frac 1N(3g-3)\right],&\text{$g\ne1$ and  $N>1$\ ,}
\endcases
\\
&
\\
K&=\cases g,&\text{$g\ne 1$, $N=1$}
\\
 2+\left[\frac 1N(g-2)\right],&\text{$g\ne1$ and  $N>1$\ .}
\endcases
\endaligned\Eq(eKL)
$$
\noindent
Here $[x]$ denotes the largest integer $\le x$.
Clearly, the first alternatives are special cases of the second ones.
Strictly speaking, the above value of $K$ is the value for
generic $n$ and $m$.
For an overall bound it has to be increased by 1 or 2 depending on $g$ and
$N$.

The algebra  $\A$ can be decomposed (as vector space) as
$$
\gathered
\A=\A_-\oplus\A_{(0)}\oplus\A_+,
\\
\A_-:=\langle A_{n,p}\mid n\le -K-1,\pN\rangle\ ,\quad
\A_+:=\langle A_{n,p}\mid n\ge 1,\pN\rangle,\quad
\\
\A_{(0)}:=\langle A_{n,p}\mid -K\le  n\le 0, \pN\rangle \ .
\endgathered\Eq(edeca)$$
and the Lie algebra $\L$ as
$$
\gathered
\L=\L_-\oplus\L_{(0)}\oplus\L_+,
\\
\L_-:=\langle e_{n,p}\mid n\le -L-1,\pN\rangle\ ,\quad
\L_+:=\langle e_{n,p}\mid n\ge 1,\pN\rangle,\quad
\\
\L_{(0)}:=\langle e_{n,p}\mid -L\le  n\le 0, \pN\rangle \ .
\endgathered\Eq(edecv)$$
Due to the almost-grading the subspaces $\A_\pm$ and $\L_\pm$
are subalgebras but
the subspaces $\A_{(0)}$, and $\L_{(0)}$ in general are not.
We use the term {\it critical strip} for the latter.

Note that $\A_+$, resp\. $\L_+$ can be described as the algebra of
functions, (resp. vector fields) vanishing at least of order one
(resp. of order 2) at the points $P_i,\iN$.
These algebras can be enlarged by allowing all elements which are
regular at all $P_i$.
This could be achieved by adding  $\ \{A_{0,p},\pN\}\ $, (resp.
$\ \{e_{0,p},e_{-1,p},\iN\}\ $ to the
set of basis elements.
We denote the enlarged algebras by
$\A_+^*$, resp. by $\L_+^*$.

On the other
hand
$\A_-$ and $\L_-$ could also be
enlarged by considering all elements which are vanishing of order one
(resp. order two) at the point $\Pif$.
Clearly they could further  be enlarged
to contain all elements which are regular at $\Pif$.
By this process  we  include elements from the critical
strip into  these algebras.

In view of Section \kzgen\ we will consider this for $\L_-$ in more
detail.
By \equ(ordi) we see that passing from $f_{n,p}^\lambda$ to
 $f_{n-1,p}^\lambda$ the order at $\Pif$ will increase by $N$.
By direct calculation we obtain for the vector fields
$$
-N<\ord_{\Pif}(e_{-L+1,p})\le 0, \quad
0<\ord_{\Pif}(e_{-L,p})\le N,\quad
N<\ord_{\Pif}(e_{-L-1,p})\le 2N\ .
$$
In particular the elements in $\L_-$ are all vanishing of at least
second order at $\Pif$.
Note that on every degree the order of $e_{n,p}$
at $\Pif$ does not depend on $p$.
Fix such a degree $n$. Due to the generic choice of the points
a suitable linear combination of $e_{n,p}$ and
$e_{n,s}, s\ne p$ has exact order 1 more than $e_{n,p}$ at $\Pif$.
In this way it is possible to find basis elements
$g_{n,i}$ of  $\Fln$ with
$$
\ord_{\Pif}(g_{n,i})=\ord_{\Pif}(e_{n,p})+i,\quad i=0,\ldots,N-1\ .
$$
After making a change of basis in this sense we will find
basis elements which can be added to the generators of
$\L_-$ and obtain a bigger algebra $\L_-'$ containing all vector
fields with vanishing order at least two at $\Pif$.
This involves only the elements of the critical strip
of
degree lower than $n$.
The corresponding basis elements have to be removed from the
critical strip.
The remaining subspace we will call  {\it reduced critical strip}
 $\L_{(0)}'$.
Clearly $\L_-'$ can be extended to $\L_-''$ containing all vector
fields vanishing at least of order one at $\Pif$  and
$\L_-^*$  containing all vector fields regular at $\Pif$.
The dimension of the corresponding critical strips
are one (resp. two) less than the dimension of the reduced
critical strip.
By counting the orders of the basis elements its dimension can be calculated
as
$$
\dim \L_{0}'=N+N+(3g-3)+1+1=2N+3g-1\ .\Eq(dimcs)
$$
The first two terms correspond to $\L_0$ and $\L_{-1}$. For
$g\ge 2$ the intermediate term comes from the basis
vector fields which have poles
at the $P_i,\iN$ and $\Pif$.
The $1+1$ corresponds to  the above constructed basis vector fields
with exact order one, resp. two at $\Pif$.
In Section \kzgen\ we will explain how the elements of the
critical strip (resp. subsets of them) are related to
tangent directions in the moduli space.

A similar decomposition is valid for the critical strip
of the function algebra. This yields a modified
$\A_-'$ and $\A_{(0)}'$. Note that
$$
N\le\ord_{\Pif}(A_{-K-1,p})< 2N\ .
$$
Let us call the subalgebra of functions which are regular at $\Pif$
by $\A_-^*$.

\subhead
(c) Central extensions and affine  algebras of higher genus
\endsubhead

The function algebra $\A$ (considered as an abelian Lie algebra)
can be centrally extended to a Lie algebra $\Ah$ via the
Lie algebra cohomology
2-cocycle
$$
\gamma(g,h):=\cintt gdh\ .\Eq(ecfac)
$$
More precisely, $\Ah=\C\oplus\A$ as vector space  with
Lie algebra structure given by
$$
[\hat g,\hat h]=\gamma(g,h)\,t_1,\quad [t_1,\Ah\,]=0\ ,\Eq(ecfa)
$$ where we used the notation
$\ \hat g:=(0,g),\ \hat h:=(0,h),\ t_1=(1,0)$.

To obtain central extensions of the vector field algebra
(generalizing the Virasoro central extension) we have first
to choose a global holomorphic projective connection  $R$.
The defining 2-cocycle is given as
$$
\chi_R(e,f):=\cinttt\left(\frac 12(e'''f-ef''')
-R\cdot(e'f-ef')\right)dz\ .\Eq(ecva)
$$
This cocycle was introduced for the $(1,1)$ case by Krichever and
Novikov. As shown in \cite{\rSDiss} it can be extended to the
multi-point situation.
It defines a central extension $\Lh=\Lh_R$.
Another choice of the projective connection
(even if we allow meromorphic projective connections with
poles only at the points in $A$)
yields a cohomologous cocycle, hence an equivalent central
extension. We denote the non-trivial central generator by $t_2$.
The above cocycles fulfil the important following {\it locality conditions}.
\proclaim{Proposition  \kznum.5}
\cite{\rSDiss}
There are constants $T$ and $S$ such that for all $m,n\in\Z$
$$
\align
\gamma(A_{n,r},A_{m,p})\ne 0 &\implies  T\le |m+n|\le 0,
\\
\chi_R(e_{n,r},e_{m,p})\ne 0 &\implies  S\le |m+n|\le 0,
\endalign
$$
\endproclaim
By considering the order of the integrands in \equ(ecfac) and \equ(ecva)
we see that the cocycles restricted to the subalgebras
$\A_+,\A_-,\A_-'$, resp. $\L_+,\L_-,\L_-'$ are vanishing.

Again explicit expressions for $T$ and $S$ can be given, but are not
of interest here.

By the locality of
the cocycle and by setting $\ \deg(t_1):=\deg(t_2):=0\ $
we obtain in this way an almost grading for
$\Lh$ and $\Ah$.
By the vanishing of the cocycles on the subalgebras $\A_{\pm}$ and $\L_{\pm}$
they
can be identified in a natural way with
subalgebras $\Ah_{\pm}$ and $\Lh_{\pm}$ of $\Ah$, resp\. $\Lh$.
\bigskip
Let $\g$ be a reductive finite-dimensional Lie algebra with a fixed
invariant, nondegenerate symmetric bilinear form $(..|..)$,
e.g\. for the semi-simple case the Cartan-Killing form.
The {\it higher genus loop algebra} or  {\it higher genus
current algebra} is defined as
$$
\gb:=g\otimes \A,\quad
\text{with Lie product}\quad
[x\otimes g,y\otimes h]:=[x,y]\otimes g\cdot h\ .
$$
It has a central extension  $\gh:=\C\oplus \gb$ with
Lie product
$$[\widehat{x\otimes f},\widehat{y\otimes g}]=
\widehat{[x,y]\otimes (f g)}-(x|y)\cdot\gamma (f,g)\cdot t,\qquad
[\,t,\Gh]=0\ ,\Eq(eaff)
$$
(where we set $\widehat{x\otimes f}:=(0,x\otimes f)$).
For the proofs, see \cite{\rSSS}.
This algebra is called the higher genus (multi-point)
affine Lie algebra (or \KN\ algebra of affine type).
Again, we can define an almost grading  on
$\gb$ and $\gh$ by setting
$$
\deg(t):=0,\quad
\deg(\widehat{x\otimes A_{n,p}}):=\deg(x\otimes A_{n,p}):=n
$$
and obtain a splitting as above
$$
\gathered
\G=\G_-\oplus\G_{(0)}\oplus\G_+,\quad\text{with}\quad
\G_\beta=\g\otimes \A_\beta,\quad \beta\in\{-,(0),+\}\ ,
\\
\Gh=\Gh_-\oplus\Gh_{(0)}\oplus\Gh_+\quad\text{with}\quad
\Gh_\pm\cong\G_\pm\quad\text{und}\quad\Gh_{(0)}=\G_{(0)}\oplus\C\cdot t\ .
\endgathered
\Eq(espaff)
$$
In particular, $\Gh_\pm$ and $\G_\pm$ are subalgebras.
The corresponding is true for the enlarged algebras.
Of special interest are
$$
\gh_-^*=\gb_-^*=\g\otimes \A_-^*,\qquad
\gh_+^*=\gb_+^*\oplus \C t=(\g\otimes \A_+^*)\oplus\C t\ .
$$

\proclaim{Lemma \kznum.6}
$$
1=\sum_{p=1}^N A_{0,p}\ .\Eq(one)
$$
\endproclaim
\demo{Proof}
Using \equ(edu) we can write
$$
1=\sum_{n\in\Z}\,\sum_{p=1}^N \kndual {1}{\w^{n,p}} A_{n,p}\ .
$$
Calculating the orders of the integrand we obtain
$$
\kndual {1}{\w^{n,p}}=\cintt \w^{n,p}=0\ ,\quad \text{for } n\ne 0\ ,
\qquad
\kndual {1}{\w^{0,p}}=1\ .
$$
The latter relation is due to the normalization.
Hence the claim.\qed
\enddemo
By this we see that the finite-dimensional Lie algebra $\g$ can  naturally be
considered as subalgebra of $\Gb$ and $\Gh$. It lies in the
subspace $\G_0$.

\vskip 1cm
\head
  \kzrep . Representations of the multi-point \KN\ algebras
\endhead
\vskip 0.2cm
\numsec=\kzrep   
\numfor=1          

\def\kznum{\kzrep}
Let us start this section with a definition of what we mean under
Verma modules of an affine multi-point \KN\ algebra $\gh$. The construction
presented below is a generalization of that proposed in \cite\rShns .

Let $\g$ be a simple finite-dimensional Lie algebra.
Let $\gb_0\subset\gb$ be the linear subspace of elements of degree
$0$, i.e.  $\gb_0= \bigoplus\limits_{p=1}^N \g\otimes A_{0,p}$.
Let $\gh_+=\gb_+\subset\gb$ be the linear subspace of elements of
a positive degree (see the notation \equ(espaff)).
Let $Z$ be the one-dimensional subspace of $\gh$ generated by the
central element $t$. The degree of $t$ was defined to be $0$,
hence $\gh_0=\gb_0\oplus Z$.
Recall that from Lemma \kzkn.6 it follows that  $\g\subset\gh_0$
(as a Lie subalgebra).
\proclaim{Lemma \kzrep .1} The direct sum $\gh_+':=\gh_0\oplus Z\oplus\gh_+$
 is a Lie subalgebra of $\ \gh$.
\endproclaim
\demo{Proof}
Recall from Section \kzkn\ that $\A_0\oplus\A_+$ is the (associative)
subalgebra of functions which are regular at all the $P_i$, $\iN$.
Hence the claim follows. See also Proposition  \kzkn.4.
\qed
\enddemo

Consider the direct sum $\g_{(N)}$ of $N$ copies of the Lie algebra
$\g$: $\g_{(N)}:= \g_1\oplus\ldots\oplus\g_N$ where an
isomorphism $\vp_p:\g_p\rightarrow\g$ is fixed
for each $p=1,\ldots,N$. In what
follows we shall assume that the maps $\vp_p$ have a special
structure. Namely, let $G:=\exp\g$ be the associated Lie group. Set
 $\gamma:=(\gamma_1,\ldots,\gamma_N )$ with $\gamma_p\in G$, $p=1,\ldots,N$
arbitrary and take $\vp_p:=Ad\,\gamma_p$.

Choose and fix a Cartan subalgebra $\car\subset\g$, a corresponding
Borel subalgebra $\bor\subset\g$ and a corresponding upper nilpotent
subalgebra $\nil\subset\g$.  Set $\car_p:= \vp_p^{-1}\car$,
$\bor_p:= \vp_p^{-1}\bor$, $\nil_p:= \vp_p^{-1}\nil$ and
$\bor_{(N)}:= \bor_1\oplus\ldots\oplus\bor_N$. Let $V_p$ be a
one-dimensional linear space over $\Bbb C$ with a fixed basis vector
$v_p$ ($p=1,\ldots,N$) and set $V:=\bigotimes\limits_{p=1}^N V_p$. Choose
an arbitrary $N$-tuple $\l:=\{\l_1,\ldots,\l_N\}$ with
$\l_p\in\car_p^*$, $p=1,\ldots,N$.

Define for  $p=1,\ldots,N$ the one-dimensional
 representations of $\bor_p$ on $V_p$
by
  $$ h_pv_p=\l_p(h_p)v_p,\quad n_pv_p=0,\qquad
  \text{for}\quad h_p\in\car_p,\quad n_p\in\nil_p.        \teq(verma)
  $$
Extend this to an
one-dimensional representation of the Lie algebra $\bor_{(N)}$ in the
linear space $V$ by decomposing
$x_{(N)}:=x_1\oplus\ldots\oplus x_N$  with
$x_1\in\bor_1,\ldots,x_N\in\bor_N$ and setting
  $$ \align
     x_{(N)} (v_1\otimes\ldots\otimes v_N) =
     &(x_1v_1)\otimes v_2\otimes\ldots\otimes v_N +\teq(tenz)\\
     &v_1\otimes (x_2v_2)\otimes\ldots\otimes v_N +
     \ldots +
     v_1\otimes v_2\otimes\ldots\otimes (x_Nv_N).
     \endalign
  $$
Let us denote this representation of the Lie algebra $\bor_{(N)}$ in
the linear space $V$ by $\tau_{\l,\gamma}$.

Set $\bb_0=\bigoplus\limits_{p=1}^N \bor\otimes A_{0,p}\subseteq \gb_0$.
Clearly, $\bh :=\bb_0\oplus Z\oplus\gh_+$ is a Lie subalgebra of $\gh$.
Take $\l$ and $\gamma$  (arbitrary) $N$-tuples as  defined above
and choose $\delta\in\C$.
Our next goal is to assign  to each triple $(\l,\gamma,\d)$ a one-dimensional
representation of the Lie algebra $\bh$.
First we define a linear map
$$
\gathered
\vp :\bh\rightarrow \bor_{(N)},\quad\text{with }
\vp_{|Z}=0,\quad \vp_{|\gh_+}=0 \quad\text{and},\\
x=
\bigoplus\limits_{p=1}^N x_p\otimes A_{0,p}\ \in\bb_0\quad\mapsto
\vp^{-1}_1(x_1)\oplus\ldots\oplus \vp^{-1}_N(x_N)
\endgathered
$$
\proclaim{Lemma \kzrep .2} The map $\vp$ is a Lie homomorphism.
\endproclaim
\demo{Proof}
As special cases of  the relations in Proposition  \kzkn.4 we obtain
  $$\align A_{0,p}A_{0,p} &= A_{0,p} + B_p,\qquad   B_p\in  \A_+\teq(ApAp)\\
           A_{0,p}A_{0,q} &= B_{p,q}  \qquad\qquad
  B_{p,q}\in\A_+,\quad q\ne p\ .
   \teq(ApAq)
    \endalign
  $$
The cocycle for defining the central extension
vanishes on $\A_0\oplus\A_+$.
This implies the claim.
\qed
\enddemo
Let us note that the map can be extended to $\ \gh_0\oplus\gh_+\to
\g_{(N)}\ $
by the same definition.

On the one-dimensional space $V$ we define the representation
$\tau_{\l,\gamma,\d}$ of $\bh$  by setting
$$\tau_{\l,\gamma,\d}(x_0\oplus x_+\oplus t):=
\tau_{\l,\gamma}(\varphi(x_0))+\delta\cdot Id\ ,
$$
with respect to the decomposition of $\bh$.
\definition{Definition \kznum .3} The linear space
 $$
\widehat{V}_{\l,\gamma,\d}:=U(\gh)\otimes_{U(\bh)}V\Eq(vermam)
$$
 with its natural structure of a
 $\gh$-module is called the Verma module of the Lie algebra $\gh$
 corresponding to the data $(\l,\gamma,\d)$. As usual
 $U(\cdot)$ denotes the universal enveloping algebra of the
 corresponding Lie algebra. $U(\bh)$ operates on $V$ via the
 representation $\tau_{\l,\gamma,\d}$. The $N$-tuple $\l$ is called
 the {\it weight} of the Verma module, the elements $\gamma_p\in\gamma$
 ($p=1,\ldots,N$) are called the {\it twists} and the complex number $\d$
 is called the {\it level} of the Verma module.
The vector $v_{\l,\gamma,\d}=v_1\otimes\cdots\otimes v_N\ $ is
called the {\it highest weight vector} of the module $V_{\l,\gamma,\d}$.
\enddefinition
\proclaim{Proposition  \kznum.4}
The $\gh$-module  $\widehat{V}_{\l,\gamma,\d}$ is a $\g$-module and
contains the
module
$$
V_{\l,\gamma}:= V_{\l_1}\otimes V_{\l_2}\otimes\cdots\otimes
V_{\l_N}\ ,
$$
where  for each $p=1,\ldots,N$, $V_{\l_p}$ is the highest weight
module of $\g_p$ of weight $\l_p$ and $\g$ operates on $V_{\l_p}$
when twisted by the automorphism $Ad\,\gamma_p$.
\endproclaim
\demo{Proof}
Recall that $\frak g$ is embedded as a subalgebra into $\gh$ via
$x\mapsto x\otimes 1=\sum_{p=1}^N x\otimes A_{0,p}$.
Hence,  the module  $\widehat{V}_{\l,\gamma,\d}$ is also a
 $\frak g$-module. The $\g$ submodule generated from
the highest weight vector $v_{\l,\gamma,\d}$ is the module
$V_{\l,\gamma}$. From the definition of the $\g$-action
on the $V_{\l,\gamma,\d}$ follows that the
action on $V_{\l,\gamma}$
for decomposable tensors
$w=w_1\otimes w_2\cdots\otimes w_n$
is given as
$$
g.w:=(\vp_1(g).w_1)\otimes w_2\cdots\otimes  w_N+\cdots+
w_1\otimes\cdots \otimes(\vp_N(g).w_N)\ .\qed
$$
\enddemo

We have inside the representation space further important
subspaces.
Recall that $\gh_+'$ is a Lie subalgebra (see Lemma \kznum.1)
of $\gh$. The subspace of $\widehat{V}_{\l,\gamma,\d}$
generated as $\gh_+'$-submodule from   $v_{\l,\gamma,\d}$
is called the {\it subspace of degree zero}.
It is denoted by  $\widehat{V}_{(\l,\gamma,\d),0}$.
Clearly,  $V_{\l,\gamma}\subseteq \widehat{V}_{(\l,\gamma,\d),0}$
Note also that $\gh_+$ annihilates $\widehat{V}_{(\l,\gamma,\d),0}$.

$\gh$ contains also the subalgebra $\gb_{-}^*\cong\gh_{-}^*$.
Note that this subalgebra was defined in Section \kzkn\ using $x\otimes f$
with $f\in\A$  and $f$ regular at $\Pif$.
The same is true for  $\gb_{-}'\cong\gh_{-}'$
where we require that the functions have zeros at $\Pif$.
It is possible to define {\it conformal blocks} as
the space of coinvariants
$$
\widehat{V}_{\l,\gamma,\d}/\gb_{-}^*\widehat{V}_{\l,\gamma,\d}\ .
\Eq(cbdef)
$$
We will postpone the discussion of the conformal blocks,
their structure,
etc. to the forthcoming part II of this article.
Note that it is also possible to replace in
\equ(cbdef)  $\gb_{-}^*$ by  $\gb_{-}'$ and a bigger space will be obtained.

\bigskip
Let us now consider more general modules over $\gh$.
\definition{Definition \kznum .5}
 A module $\Vh$ over the Lie algebra $\Gh$ (resp\. a
 representation of  $\Gh$) is called an {\it admissible module} (resp\.
an {\it
 admissible representation})
if the central element $t$ operates as $c\cdot Id$ with $c\in\C$ and
if for every $v\in {\Vh}$ and for all $x\in\g$
 one has $x(n)v=0$ for $n\gg 0$.
\enddefinition
It is evident that each Verma module \equ(vermam)
is an admissible module.

For each admissible representation of the affine \KN\ algebra $\gh$
the (affine) Sugawara construction yields a
representation of the centrally extended \KN\ vector field
algebra $\Lh$ (the Virasoro type algebra).
This representation is called {\it Sugawara   representation}.
 The abelian version for the two-point case
was introduced in \cite{\rKNFb}. The
nonabelian case
was later considered in \cite{\rBono}, \cite{\rSSS}.
In \cite{\rSSS} also the multi-point version was given.

Let $\Vh$ be a fixed admissible module.  For
$u\otimes A_{n,p}$
with $A_{n,p}$ being a basis element of the algebra $\Ah$
and $u\in\g$
we will denote the corresponding
operator in $\Vh$ by $\ u(n,p)\ $ as well as by $u(A_{n,p})$.

Recall that we assume $\g$ to be a finite dimensional simple Lie
algebra. We choose a basis $\ u_i,\ i=1,\ldots,\dim\g\ $ of $\g$ and
the corresponding dual basis $\ u^i,\  i=1,\ldots,\dim\g$ with
respect to the invariant non-degenerate  symmetric bilinear form
$(..|..)$.  The Casimir  element $\ \Omega^0=\sum_{i=1}^{\dim \g}
u_iu^i\ $ of the universal enveloping algebra  $U(\g)$ is independent
of the choice of the basis.
In an abuse of notation we denote  $\sum_i u_i(n,p)u^i(m,q)$ simply
by $u(n,p)u(m,q)$.

We define the higher genus {\it Sugawara operator}
(also called {\it Segal operator} or {\it energy-momentum tensor operator})
as
$$T(Q):=\frac 12\sum_{n,m}\sum_{p,s}
\nord{u(n,p)u(m,s)}\w^{n,p}(Q)\w^{m,s}(Q)\ .\Eq(suga)
$$
By $\ \nord{....}\ $ we denote some normal ordering.
The summation  here and in the following
formulas  for the first indices $n,m$ are over
$\Z$ and for the second indices $p,s$ over $\{1,\ldots,N\}$.
The precise form of the normal ordering is not of importance.
As an example we may take the following ``standard
normal ordering''  ($x,y\in\g$)
$$
\nord{x(n,p)y(m,r)}\ :=\cases x(n,p)y(m,r),&n\le m
                                           \\
                            y(m,r)x(n,p),&n>m\ .
                      \endcases
\Eq(normst)
$$
The expression $T(Q)$
can be considered as formal series
of quadratic differentials in the variable $Q$ with
operator-valued coefficients.
Expanding it over
the basis $\Omega^{k,r}$ of the
quadratic differentials we obtain
$$
  T(Q)=\sum_k\sum_r L_{k,r}\cdot\Omega^{k,r}(Q)\ ,\Eq(sugb)
$$
with
$$
  \gathered
  L_{k,r}=\cintt T(Q)e_{k,r}(Q)=\frac 12\sum_{n,m}\sum_{p,s}
  \nord{u(n,p)u(m,s)}l_{(k,r)}^{(n,p)(m,s)},\\
  \text{where}\qquad
  \ l_{(k,r)}^{(n,p)(m,s)}:=\cintt \w^{n,p}(Q)\w^{m,s}(Q)e_{k,r}(Q)\ .
  \endgathered\Eq(sugc)
$$

The following theorem was proved in \cite\rSSS .
\proclaim{Theorem \kznum.6}
Let $\g$
be a finite dimensional either abelian or simple Lie algebra and
$2\k$ be the eigenvalue of its Casimir operator in the adjoint
representation.
Let $V$ be  an admissible representation  where the central element
$t$ operates as  $\ce\cdot identity$. If $\ce+\k\ne 0$ then the rescaled
modes
$$L_{k,r}^*=\frac {-1}{2(\ce+\k)}\sum_{n,m}\sum_{p,s}
 \nord{u(n,p)u(m,s)}l_{(k,r)}^{(n,p)(m,s)}
\ ,\Eq(sugm)$$
of the Sugawara operator
are well-defined operators on $V$ and define a representation of
the centrally extended vector field algebra $\Lh$
(the Virasoro-type algebra).
\endproclaim
We
call the $L_{k,r}^*$, resp. the $L_{k,r}$ the Sugawara
operators too.  The representation obtained in this way is the
Sugawara representation of the Lie algebra $\Lh$ corresponding to the
given admissible representation $\Vh$ of $\gh$.

\vskip 1cm
\head
  \kzgen.     Moduli of curves with marked points and
the general form of the KZ equation
\endhead
\vskip 0.2cm
\numsec=\kzgen   

\def\kznum{\kzgen}
\numfor=1      
%
\subhead
(a) Moduli spaces
\endsubhead

Let $\MgN$, resp\. $\MgNe$ be the moduli space of smooth, projective
curves of genus $g$ (over $\C$)
with $N$ (resp. $N+1$) marked  ordered distinct points.
Equivalently, it
can be described as moduli space of compact
Riemann surfaces with marked points.
A point in $\MgNe$ is given by the equivalence classes of
the data $\mpt$ with $M$ a smooth, projective curve and $P_i\in M,\iN,\infty$.
Two such tuples are identified if they are  isomorphic
(under algebraic  maps) as  curves with marked points.
Let us denote the equivalence class by $[..]$.
To avoid special considerations we will mainly assume in this
section $g\ge 2$.
In this article we are only dealing with the local situation at
a generic curve $M$ with generic markings $(P_1,P_2,\ldots, P_N,\Pif)$
\footnote
{
Sometimes one understands by a marking of a curve the choice of a
symplectic basis for the homology. But in this article a marking
refers always to the choice of (ordered) points on the curve.
}.
Hence, it is enough to consider a small open subset $\Wt$ around the
point $\tilde b=[\mpt]$.
A generic curve of $g\ge 2$ does not admit nontrivial
infinitesimal automorphisms
and  we may assume that there exists over $\Wt$ a universal
family of curves with  marked points.
In particular, there is a proper, flat family of smooth curves over  $\Wt$
$$\Cal U\to  \Wt\ ,\Eq(unim)$$
such that
for the points $\tilde b=[\mpt]\in\Wt$
we have $\ \pi^{-1}(\tilde b)= M$ and that the sections
defined as
$$
\sigma_i:\Wt\to\Cal U,\quad\sigma_i(\tilde b)=P_i,\quad \iN,\infty
\Eq(sect)
$$
are holomorphic.
By ``forgetting the marking'' we obtain a (local) map
 $\tilde\nu:\Wt\subseteq \MgNe\to \Mgo$.
Note that the family \equ(unim) is the pullback of the universal
family over $\Mgo$ under $\tilde\nu$.
For more background information, see \cite{\rUcft, Sect. 1.2, Sect. 1.3},
in particular see Thm. 1.2.9 of \cite{\rUcft}.

Let us fix a section $\sinf$
of the universal family of isomorphy classes of curves (without marking).
In particular, for every curve there is a point chosen in a
manner depending smoothly on the moduli.
(Recall, we are only dealing with the local and generic situation.)
The analytic subset
$$
W':=\{\tilde b=[\mpt]\mid \Pif=\sinf([M])\}\quad\subseteq\quad \Wt\Eq(msidw)
$$
can be identified with an open subset $W$ of $\MgN$ via
$$
\tilde b=[ (M,P_1,P_2,\ldots, P_N,\sinf([M]))]
\ \to\
b=[ (M,P_1,P_2,\ldots, P_N)]\ .\Eq(msid)
$$
(By the genericity the map is 1:1.)

The dimensions of the moduli spaces are well-known
(see also further down in this section)
$$
\dim_b(\MgN)=\cases      3g-3+N,& g\ge 2\\
                          \max(1,N),& g=1\\
                          \max{(0,N-3)},&g=0\ .
             \endcases
\Eq(dimM)
$$
Note that for $N\ge 3-2g$ the first expression is valid for every genus.

Over the points $b\in W$ we can apply the construction
of Section \kzkn\ and Section \kzrep\ and obtain objects
over the subset $W$ of $\MgN$
$$
\A_{b},\ \Ah_{b},\ \L_{b},\ \Lh_{b},\ \gh_{b},\ \Fl_{b}\ ,\Eq(locob)
$$
depending on the points $b\in W\subseteq\MgN$.
We need a sheaf description of some of the objects.
Consider the universal family $\pi:\U\to W$.
First let us introduce the notation
$\ S_b=\sum_{i=1}^N P_i\ $ for the divisor on $M$ corresponding to
the moduli point $b=[\mpp]$.
Varying $b$ this defines
the divisor of sections
$\ S=\sum_{i=1}^N \sigma_i(W)\ $
in the family $\Cal U$.
We have to enlarge the divisors by adding the reference point
$\Pif=\sinf([M])$, resp. the section $\sinf$.
Denote by $\nu:W\subseteq\MgN\to\Mgo$ the map obtained by forgetting
the marking. We set
$$
\widetilde{S}_b=\sum_{i=1}^N P_i+\sinf(M),\quad
\widetilde{S}=\sum_{i=1}^N \sigma_i(W)+\sinf(\nu(W))\ .
\Eq(secdive)
$$
Denote by $\O_{\U}$ the sheaf of regular functions on $\U$.
As usual  set $\O_{\U}(k\St)$, $k\in\Z$ the sheaf of functions which have
zeros of order at least $-k$ along the
divisor $\St$. In particular, for $k\in\N$ this says
that the functions have poles of order at most $k$ at the divisor
$\St$.
By $\O_{\U}(*\St)$ we understand the sheaf of functions which
have poles along the divisor $\St$.
Also we set
$$
\pi_*\O_{\U}(*\St):=\lim_{k\to\infty}\pi_*(\O_{\U}(k\St))\ .
\Eq(piste)
$$
It is a locally free $\O_W$-sheaf
(with $\O_W$ the structure sheaf of the space $W$).
Its vector space fibre over $b\in W$ can be identified with
$$
\Ho(M_b,\O_{M_b}(*\St))=\A_b\ ,\quad M_b=\pi^{-1}(b)\ .
$$
The sheaf $\pi_*\O_{\U}(*\St)$ can be made to a sheaf of
commutative associative $\O_W$-algebras by fibre-wise
multiplication.
We
denote
it $\A_W$.
As a sheaf of abelian Lie algebras it can be centrally extended to the
$\O_W$-sheaf
$$
\Ah_W:=\widehat{\pi_*\O_{\U}(*\St)}:=\pi_*\O_{\U}(*\St)\oplus \O_W\cdot t,
$$
where $t$ is the central element and the structure is defined as follows.
The elements $f,g\in\A_W(U)$  can be represented as functions
on
$\pi^{-1}(U)$ with  poles only along $\St$.
Let $\gamma$  be the cocycle \equ(ecfac).
Recall that the latter can be given by calculating residues
along $S$.
Then with $r,s\in\O_W(U)$
$$
[f+r\cdot t,g+s\cdot t]:=\gamma(f,g)\cdot t,
$$
defines an element of $\Ah_W(U)$.
Note that $\gamma(f,g)\in\O_W(U)$.

This construction can be extended to the affine algebra situation.
\definition{Definition \kznum.1}
Given a finite-dimensional Lie algebra $\g$
the {\it sheaf of the associated loop algebra} (or current algebra)
$\Gb_W$ and the {\it sheaf of the associated affine algebra} $\Gh_W$
are defined as
$$
\Gb_W:=\A_W\otimes\g,\qquad
\Gh_W:=\Gb_W\oplus\O_W\cdot t\ ,
\Eq(shaff)
$$
where the Lie structure is given by the naturally extended
form of \equ(eaff) (resp. without its central term for $\Gb_W$).
\enddefinition
Clearly, these are $\O_W$-sheaves.
Let $b\in W$ and let  $\O_{W,b}$ be the local ring at $b$
and $M_b$ its maximal ideal. Set $\C_b\cong  \O_{W,b}/M_b$,
then we obtain the following canonical isomorphisms
$$
\C_b\otimes\A_W\cong \A_b,\qquad
\C_b\otimes\Ah_W\cong \Ah_b,\qquad
\C_b\otimes\Gb_W\cong \Gb_b,\qquad
\C_b\otimes\Gh_W\cong \Gh_b \ .
$$
\definition{Definition \kznum.2}
A sheaf $\frak V$ of $\O_W$-modules
is called a {\it sheaf of representations} for the affine
algebra $\Gh_W$ if
the $\frak V(U)$ are
modules over $\Gh_W(U)$.
\enddefinition
For a sheaf of representations $\frak V$ we obtain that
$\frak V_b$ is a module over $\Gh_{b}$
for every point in $W$.
\medskip

We have to increase the moduli spaces by considering also first order
infinitesimal neighbourhoods around the points $P_i,\iN,\infty$.
We denote this moduli space by $\MgNn$, resp.  $\MgNen$.
The elements of  $\MgNen$ are given as
$$
\btn=[\mptn]\ ,
$$
where for $\iNi$ the additional data $z_i$ is a coordinate at $P_i$ with
$z_i(P_i)=0$. Two such tuples $\btn$ and ${{\btn}}{}^{'}$ are identified if
they are equivalent as truncated elements in $\MgNe$
and (after this identification) we have
$$
z_i'=z_i+O(z_i^2), \qquad \iNi\ .
$$
The additional degrees of freedom in the moduli space
at a point $\bt$ is given by multiplying
fixed coordinates
in $\bt$ at the $P_i'$s by
non-zero constants.
Denote the corresponding space lying above $\Wt$ by $\Wtn$.

Again, after fixing a first order infinitesimal neighbourhood
around the section $\sinf$, we can identify the subspace defined
similarly to \equ(msidw) with an open subspace of
$\MgNn$ containing a neighbourhood of the point at which we make our
consideration. Denote this subspace by $\Wn$.
Clearly we have $N$ degrees of additional freedom and hence
 it follows from \equ(dimM)
$$
\dim_{b^{(1)}}(\MgNn)=\cases      3g-3+2N,& g\ge 2\\
                          2N,& g=1\\
                          \max{(0,2N-3)},&g=0\ .
             \endcases
\Eq(dimMn)
$$
The map $\ \eta:\Wn\to W\ $ obtained by forgetting the coordinates
is a surjective analytic map.
Hence by pulling back via $\eta$ the  objects \equ(locob)
and the sheaves
$\ \A_W, \Ah_W, \Gb_W, \Gh_W $ we obtain
sheaves
$\ \A_{\Wn}, \Ah_{\Wn}, \Gb_{\Wn}, \Gh_{\Wn} $
over $\Wn$.
Moreover, over $\Wn$ also the basis elements $f_{n,p}^\l$ are well-defined.
This was the reason for enlarging the moduli space.
Recall that for fixing the basis elements
 a choice of coordinates $z_i$ around the points $P_i$
were necessary.
But note also that only the class of $z_i$ under the equivalence
in $\MgNn$ is of importance.
Due to the explicit description \cite{\rSLb} they depend
analytically on the moduli.
Pulling back a sheaf of representation  $\frak V$ over $W$ we obtain
a sheaf of representation  $\frak V^{(1)}=\eta^*\frak V$
of $\Gh_{\Wn}$.
More generally, we can define sheaves of representations
over $\Wn$  directly. In particular for these sheaves of representations
operators depending on the KN basis  are well-defined.

We want to study how a different choice of coordinates
can be expressed on $\MgN$.
Take $b\in\MgN$ and choose  coordinates $z_i$ at the points $P_i$.
Only the coordinate classes are of importance. Hence we can express this
as choosing a lift $W\to W^{(1)}$ of $\eta$.
{}From the construction of the KN basis elements in Section \kzkn\ the
following lemma is immediate.
\proclaim{Lemma \kznum.3}
Let $\ z_p'=\alpha_p\cdot z_p+O(z_p^2)$ be another coordinate
at $P_p$. Let $f_{n,p}^{\l}$ ($f_{n,p}^{\l\prime}$)
be a KN basis element of $\Fl$ w.r.t. $z_p$ (w.r.t. $z_p'$) then
$$
f_{n,p}^{\l\prime}=(\alpha_p)^n f_{n,p}^\l,\quad\text{and }
f_{n,s}^{\l\prime}=f_{n,s}^\l,\quad s\ne p\ .\Eq(bta)
$$
\endproclaim
Note that the $\alpha_p$ are nonvanishing local functions on $\MgN$.
This behaviour has some important consequences.

(1) The grading of $\ \gh\ $ is given with respect to the
basis elements $A_{n,p}$ in $\A$.
{}From \equ(bta) it follows that the degree is not
changed by passing from one system of coordinates to another.
Clearly, this globalizes over $W$. Hence we can equip
$\A_W, \Ah_W, \Gb_W$ and
$\Gh_W$ with an almost grading.
This allows to define a {sheaf of admissible representations}
to be a sheaf of representations where all representations are
admissible.
In addition we will usually require (if nothing else is said)
that the central element
$t$
operates as $c\cdot id$ with $c$ a  function on $W$.
This function is called the level function.
Very often we will  even assume $c$ to be a constant $\ce$, which is
just called the level of the representation.

(2) Due to the possibility to globalize the grading to
$\Gh_W$ it is possible to define the {\it Verma module
sheaf}
in a straightforward manner extending Definition \kzrep .3
$$
\widehat{V}_{(\l,\gamma,\d),W}:=U(\gh_W)\otimes_{U(\bh_W)}V\Eq(vermash)\ .
$$
On first sight it looks as if the Verma module sheaf is only defined over
$\Wn$ because the basis elements $A_{0,p}$ are involved.
But from \equ(bta) it follows that they are indeed independent
of the coordinate classes.
The Verma module sheaves are sheaves of admissible representations.
Note, it is even possible to vary the data $(\l,\gamma,\d)$
over the moduli.

The subspace of degree zero defines
a subsheaf. The same is true for the subspace which is annihilated by
$\gh_+$ and for the subspace $\gb_{-}^*\widehat{V}_{\l,\gamma,\d}$,
see \equ(cbdef).
One possible way to define the sheaf of conformal blocks
is to define the quotient sheaf with respect to
the latter subsheaf.
The discussion of these objects is postponed
to  the forthcoming part II of this article.

(3)
For every sheaf of admissible representations
over $\Wn$
the Sugawara operators are well-defined.
Note that the individual operators $u(n,p)$ (see Section
\kzrep)
depend on the coordinates.
Let us consider a sheaf of admissible representations
over $W$ and choose  coordinates. From \equ(bta) we know
how the operators transform if we choose a different coordinate
$z_p'=\alpha_p z_p+O(z_p^2)$. We obtain
$\ (u(n,p)){}'=(\alpha_p)^n u(n,p)$.
The factor will be cancelled by the contribution from
$\w^{n,p}{}'=(\alpha_p^){-n}\w^{n,p}$.
Hence it follows from \equ(suga) for the Sugawara operators
that $T(Q)'=T(Q)$.
The individual operators $L_{k,s}$ depend indeed on the coordinates.
They are only well-defined over $\Wn$.
They transform as the vector fields  $e_{k,s}$ do,
as can be seen from \equ(sugb), \equ(sugc).
If we assume  the level function  $c$ to
obey the condition $(c+\ka)\ne 0$ (see Theorem~\kzrep.6)
and if we assign to the vector field $l\in L$
the operator
$$
T[l]:=\frac {-1}{\ce+\ka}\cdot\cintt T(Q)l(Q)  \Eq(sugvf)
$$
we see that this operator does not depend on the
coordinates.  Clearly this is the Sugawara representation.
We get $T[e_{n,p}]=L^*_{n,p}$, see \equ(sugc).

\bigskip
Let us now consider the tangent space at the moduli spaces.
The Kodaira-Spencer map for a versal  family
of complex analytic manifolds $Y\to B$
over the base $B$ at the base point $b\in B$,
$$
T_b(B)\to\He(Y_b,T_{Y_b})\Eq(KSh)$$
is an isomorphism (e.g. \cite{\rKod}).
Here $T_b(B)$ denotes the tangent space of $B$ at the point $b$, $Y_b$ is
the fibre over $b$ and $T_{Y_b}$ the (holomorphic) tangent sheaf
of $Y_b$.
We are in the local
generic situation where we have a universal family.
Hence we  can employ \equ(KSh).
Let $M$ be the curve fixed by  $b$, resp $\bn$.
We obtain
$$
T_{[M]}(\Mgo)\cong \He(M,T_M)\ .
$$
More generally,
(see also \cite{\rUcft}, \cite{\rTUY})
$$
T_b(\MgN)\cong \He(M,T_M(-S_b)), \qquad
T_{\bn}(\MgNn)\cong \He(M,T_M(-2S_b)) \ .
\Eq(KShom)
$$
The first order  vanishing condition at the points $P_1,\ldots, P_N$
comes from the fact that the  vector fields which do not generate
a non-trivial  complex deformation of the curve should also
not move the points to be a trivial deformation of the
marked curve.
The second order vanishing condition corresponds to the fact that
it should additionally not change the first order
infinitesimal neighbourhood.
Via Serre duality \cite{\rHAG}, \cite{\rSRS}
we have
$$
\gathered
\He(M,T_M)\cong \Ho(M,\K^2),\qquad
\He(M,T_M(-S_b))\cong \Ho(M,\K^2(S_b)),
\\
\He(M,T_M(-2S_b)\cong \Ho(M,\K^2(2S_b))\ .
\endgathered
\Eq(serre)
$$
Note that for the Kodaira-Spencer mapping there exists also a sheaf version
\cite{\rUcft, Cor. 1.2.5}
$$
T_W\cong R^1\pi_*T_{\Cal U/W}(-S)\ .
\Eq(KSsh)
$$

\subhead
(b) \KN\ algebras and tangent vectors of the moduli spaces
\endsubhead

We want to show that the elements of the cohomology groups
in \equ(KShom) can be identified with
elements
of the critical strip
of the  \KN\ vector field algebra.
Hence the latter can be identified with tangent vectors to the moduli
spaces.

Let $M$ be the Riemann surface we are dealing with and let
$\Uif$ be a coordinate disc around $\Pif$, such that
$P_1,\dots,P_N\notin \Uif$.
Let $U_1=M\setminus\{\Pif\}$.
Because $U_1$ and $U_\infty$ are affine (resp\. Stein) \cite{\rHAG,p.297}
we get  $\He(U_j,F)=0$, $j=1,\infty$ for every coherent sheaf $F$.
Hence, the sheaf cohomology can be given  as
Cech cohomology with respect to the covering $\{U_1,\Uif\}$.
Set $\Uifs=U_1\cap\Uif=\Uif\setminus\{\Pif\}$.
The 2-cocycles can be given by $s_{1,\infty}\in F(\Uifs)$
($s_{0,0}=s_{\infty,\infty}=0$), hence by
arbitrary sections over the punctured coordinate disc $\Uifs$

Coming back to the holomorphic tangent bundle $T_M$.
For any element  $f$ of the KN vector field algebra its
restriction to $\Uifs$ is holomorphic and  defines
an element of $\He(M,T_M)$.
Note that it defines also an element of
 $\He(M,T_M(D))$,
where $D$ is any divisor supported outside of  $\Uifs$.
We introduce the map
$$
\theta_D:\L\ \to\  \He(M,T_M(D)),\qquad
f\mapsto \theta_D(f):=[f_{|\Uifs}] \ .
\Eq(theta)
$$
If the divisor $D$ is clear from the context we will suppress it in
the notation.
For us only the divisors
$\ D=-kS_b\ $ with $k\in\N_0$ are of importance.

Recall the discussion in Section \kzkn\ on the reduced critical strip.
It can be decomposed
as
$$
\L_{(0)}'=\langle e_{0,p},\pN\rangle\oplus
\langle e_{-1,p},\pN\rangle\oplus
\L_{(0)}^*\oplus \L_{(0)}^\infty\ .\Eq(decomp)
$$
Here $\L_{(0)}^*$ is the subspace generated by the basis elements
(after the explained change of basis elements)
with  poles at the $P_i$, $\iNi$ and
$\L_{(0)}^\infty$ is the two-dimensional space
generated by the basis vector fields of exact order one resp.
order zero at $\Pif$.
Note that (assuming $g\ge 2$)
$$
\dim \L_{(0)}^*=3g-3\ .
$$
Recall that $\L_-^*$ was introduced as the algebra of all vector
fields (in $\L$) regular at $\Pif$.
\proclaim{Proposition  \kznum.4}
Set $S_b=\sum_{i=1}^NP_i$. The map
$$\theta_{-kS_b}:\L\ \to\ \He(M,T_M(-kS_b),\qquad k\ge 0\Eq(crhomm)$$
is a surjective map.
It gives an isomorphism
$$
\L_{k-2}\oplus \L_{k-3}\cdots\oplus\L_{-1}\oplus
\L_{(0)}^*\ \cong\ \He(M,T_M(-kS))\ .\Eq(crhom)
$$
The kernel of the map is given as
$$
\ker \theta_{-kS_b}=\L_-^*\oplus\bigoplus_{n\ge k-1}\L_n\ .\Eq(thetaker)
$$
\endproclaim
\demo{Proof}
The map $\theta_{-kS_b}$ is linear.
Using Serre duality, we calculate
$\dim\He(M,T_M(-kS_b))=\dim\Ho(M,K^2(kS_b))$. But
$\deg(K^2(kS_b))=2(2g-2)+kN\ge 2g-1$  for $g\ge 2$.
By Riemann-Roch (in the non-special region) \cite{\rSRS}
we obtain
$$
\dim\He(M,T_M(-kS_b))=\dim\Ho(M,K^2(kS_b))=3(g-1)+kN\ .
$$
This coincides with the dimensions of the spaces
on the l.h.s. of \equ(crhom).
First we prove \equ(thetaker).
Assume $\theta_{-kS_b}(f)=0$. Hence, we can
write
$$
f_{|\Uifs}={f_1}_{|\Uifs}  - {f_0}_{|\Uifs}
$$
with $f_0$ defined and regular on $\Uif$ and
$f_1$ defined outside of $\Pif$
 and with  order at least $k$ at every point
$P_i,\iN$.
Because $f_0$ is regular at $\Pif$  the vector field
$f_1$ has to have the same
singular
part as $f$ at  $\Pif$. In particular it can be extended to a
global meromorphic vector field. This implies $f_1\in\L$.
{}From $\ord_{P_i}(f_1)\ge k$, $\iN$ it follows
$f_1\in\bigoplus_{n\ge k-1}\L_n$.
Now ${f_0}_{|\Uifs}= {({f_1}-f)}_{|\Uifs}$.
But ${f_1}-f$ is a globally defined meromorphic vector field with poles
at most at $\{P_1,\ldots, P_N,\Pif\}$.
Hence the same is true
for the extension of $f_0$
(which we denote by the same symbol).
In particular $f_0\in\L$. Due to the regularity at $\Pif$
we get $f\in\L_-^*$. This shows $\subseteq$.
\nl
For $f\in\bigoplus_{n\ge k-1}\L_n$ we set $f_1=f_{|U_1},\ f_0=0$, and
for  $f\in\L_-^*$  we set $f_1=0,\ f_0=-f_{|\Uif}$. We see that
their cohomology classes will vanish.
Hence, \equ(thetaker)
\nl
{}From \equ(thetaker) it follows that $\theta_{-kS_b}$ is
injective if restricted to the complementary space.
{}From the equality of the dimension follows   surjectivity
and \equ(crhom).
\qed
\enddemo
Note that the spaces $\ker\theta_{-kS_b}$  are invariantly defined.
The following (linear) isomorphism
are of special interest for us
$$
\gathered
\He(M,T_M)\cong \L_{(0)}^*,\qquad
\He(M,T_M(-S_b))\cong \L_{-1}\oplus\L_{(0)}^*, \qquad
\\
\He(M,T_M(-2S_b))\cong  \L_{0}\oplus\L_{-1}\oplus\L_{(0)}^*\ .
\endgathered
\Eq(crhomd)
$$
Using the Kodaira-Spencer map \equ(KShom) we obtain
\proclaim{Theorem \kznum.5}
The tangent spaces of the moduli spaces
$\Mgo$, $\MgN$ and $\MgNn$
at the points which correspond to the curve $M$ with marked points
$I=(P_1,P_2,\dots,P_N)$ and $O=(\Pif=\sinf(M))$
can be identified with  the following subspaces of the critical strip
of the Krichever-Novikov vector field algebra assigned to this marked
curve:
$$
T_{[M]}\Mgo\cong\L_{(0)}^*,\qquad
T_{b}\MgN\cong\L_{-1}\oplus\L_{(0)}^*,\qquad
T_{\bn}\MgNn\cong\L_{0}\oplus\L_{-1}\oplus\L_{(0)}^* \ ,
\Eq(tansp)
$$
where
$b=[\mpp]$ and $\bn=[\mppn]$.
\endproclaim
Because the vector fields in $\ker \theta\subset\L$ are not
corresponding to deformations in the moduli we sometimes call
this vector fields {\it vertical vector fields}.
They should not be confused with the sections of the
relative tangent sheaf of the universal family.
The following should also be kept in mind.
As already remarked, $\ker\theta$ is
invariantly defined. But the definition of the critical strip, hence
of the complementary subspace to $\ker\theta$,  is
only fixed by the order prescription
for the basis. A different
prescription (which involves changing the
required orders) will yield a different identification of
tangent vectors on the moduli space with vector fields
in the fibre.

This was a description of the connection between  the
Krichever-Novikov basis elements and the
infinitesimal  moduli parameter
using Cech cohomology.
For another (but nevertheless equivalent) approach in the $N=1$ case, see
\cite{\rGrOr}.
\medskip
Let us close this subsection in discussing the necessary
modifications for genus 0 and 1.
In these cases there are global holomorphic vector fields.
{}Hence the decomposition of the
critical strip \equ(decomp) and its identification \equ(crhomd)
are not valid anymore.
The space $\L_{(0)}^*\oplus\L_{(0)}^\infty$ (resp. for
$g=1$ a part of it) already appears as subspace of
$\L_{0}\oplus\L_{-1}$.
This is in complete conformity with the corrected dimension of
the moduli space.

Let us first consider $g=0$. It is always possible to move three distinct
points to the triple $(0,1,\infty)$ by an automorphism of $\P^1$.
If this is done there are no further  automorphisms.
Hence the moduli space $\Cal M_{0,N}$ has a non-zero dimension
exactly for $N\ge 4$. Its dimension is $\min(0,N-3)$.
It is quite useful to map the reference point $\Pif$ always to
$\infty$.
But one should pay attention to the fact that the map
\equ(msid) is not 1:1 anymore.
In other words,  replacing $W'$ by $W$ is not
possible. Note that  $\dim W'=N-2=\dim \Cal M_{0,N}+1$
for $N\ge 3$.
In this case it is even better to work with the configuration
space (and this is usually done)
$$
\widehat{W}:=\{(P_1,P_2,\ldots,P_N)\mid P_i\in\C,\ P_i\ne P_j,\
\text{for } i\ne j\}\
$$
of $N$ points and study the remaining invariance at the end.

For $g=1$ the situation is similar.
For a generic elliptic curve $E$
we always have the translations by  points of $E$
as automorphisms.
After fixing a point as the zero of the
group law on $E$ (which might be chosen as the reference point
$\Pif$) this automorphisms are not possible anymore.
The only non-trivial automorphism which remains for the generic
curve is the
involution
$x\to -x$.
Again the map \equ(msid) is not 1:1 anymore.
For example, for $N=1$ we have $\dim W'=2$ but $\dim\Cal M_{1,1}=1$.
We have to work inside $\Cal M_{1,N+1}$.
Again, it is useful to work with the ``configuration space''
picture.
Note that for higher genus in the generic situation the
moduli space locally coincides with the ``configuration space''.

\subhead
(c) The formal KZ equations
\endsubhead

Let $\frak V$ be a sheaf of admissible representations
 of the affine algebra $\Gh$
over $W$, resp. over $W^{(1)}$ as introduced above.
We assume  the level $\ce$ to be constant and
obeying the condition $(\ce+\ka)\ne 0$.
By Theorem \kzgen.5 the elements of  the (fixed) critical strip
(resp. of  a subspace of the critical strip) correspond to tangent
vectors along the moduli space
$\MgN$, resp\. along $\MgNn$.
In particular, the maps $\theta$ introduced in \equ(theta)
are  isomorphisms if restricted to the subspaces.
Denote the tangent vectors by $X_k$, $k=1,\ldots,3g-3+N$
(resp\. $k=1,\ldots,3g-3+2N$) and set
$l_k=\theta^{-1}(X_k)$ for the corresponding element of the
critical strip.
Assume $X_k$ operates linearly as operator $\partial_k$
on the space of sections.
Assume further that $\partial_k$ operates as derivation
$$
\partial_k(s\,\Phi)=(X_k(s))\Phi+s\,\partial_k(\Phi),
$$
for $\Phi$ a section of $\frak V$ and $s$ a local function on
the base $W$, resp. on $\Wn$.

Let us concentrate on sheaves of representations
which are defined over $\MgN$. Examples are the Verma module sheaves.
As seen in Section \kzrep\  the elements of (the centrally
extended) vector field algebra operate vertically on the
fibre of the representation sheaf
 via the Sugawara representation.
Recall the definition \equ(sugvf) of the the operator $T[l]$
which is defined for
every vector field $l\in\L$.
The operator does not depend on the coordinates.

We define for sections $\Phi$ of $\frak V$ and for every $k$
the operator
$$
\nabla_k\Phi:=(\partial_k+T[l_k])\,\Phi\ .
\Eq(nab)
$$
\definition{Definition \kznum.6}
The {\it formal KZ equations} are defined as
the  set of equations
$$
\gathered
\nabla_k\,\Phi=0\ ,
\\
\text{for\ }\quad k=1,\ldots, 3g-3+N\ .
\endgathered
\Eq(NablaEq)
$$
\enddefinition
\noindent
Using \equ(sugc) and
$\ l_{k}^{(n,p)(m,s)}:=\cintt \w^{n,p}\w^{m,s}l_{k}\ $
we can rewrite this as
$$ \left(\partial_k -\frac {1}{\ce+\ka} \sum_{n,m \atop p,s} l^{(n,p)(m,s)}_k
\nord{u(n,p)u(m,s)}
\right)\;\Phi =0\ ,\quad k=1,\ldots, 3g-3+N\ .
                                                \Eq(kze)
$$
In view of the decomposition \equ(tansp) the set of equations
divides up into two subsets of different interpretation.
We have $N$ equations for  moving the points.
They are related to   $e_{-1,p}$, $\pN$.
The other ones corresponding to the
 $3g-3$ elements $l_k\in\L_{(0)}^*$ (for $g\ge2$).
They are responsible for changing the complex structure of the curve.

Sometimes we will consider sections which take their values in
certain subspaces. For the Verma module sheaf
for example the subsheaf consisting
of the elements in $\frak V$ of degree zero, or the subspace of the elements
annihilated by  $\Gh_+$ (in its sheaf versions).
Or we will consider induced  actions on quotient sheaves
in such cases when the operator $\nabla_k$ maps the subspaces
which are factored out to themselves.
With this situation we will deal in the forthcoming
part II of this article.

For calculations it is often useful to express
\equ(nab), resp\. \equ(NablaEq) in terms of the KN vector field basis $e_{n,p}$
and set $X_{n,p}=\theta(e_{n,p})$.
But the definition of $e_{n,p}$ depends on the coordinates
like
$e_{n,p}'=(\alpha_p)^n e_{n,p}$. We obtain
also
$X_{n,p}'=\theta(e_{n,p}')=(\alpha_p)^n X_{n,p}$
and $T[e_{n,p}']=(\alpha_p)^n T[e_{n,p}]$.
For the operators we obtain
$\nabla_{n,p}'=(\alpha_p)^n \nabla_{n,p}$.
In particular, the set of formal
KZ equations expressed in terms of the KN basis elements will
yield an equivalent set under coordinate transformations.
Note also that the KN basis elements are fixing the tangent
directions.

The corresponding set of formal KZ equations for
sheaves of representations over $\MgNn$ are obtained by adding $N$ additional
equations corresponding to the vector fields
$e_{0,p},\ \pN$.
They deform local coordinates at the marked points.
In this case the operators $\nabla_{n,p}$ are well-defined.

Keep in mind that (for example) it is not possible to talk about changing
the complex structure and  ``fixing the points''.
Neither is it possible to move  the points
(in $\MgNn$) and ``fix the coordinates'').
Nevertheless, if we have global coordinates
we have for every movement of the points also  a
definite change of coordinates.
In other words, we are considering a special
subspace of $W^{(1)}\subseteq \MgNn$ isomorphic to $W\subseteq \MgN$.
In this case it is possible to consider the  sheaves of representation and
the corresponding Sugawara operators
which are defined over $\MgNn$ as sheaves of representations over $\MgN$
and drop the
corresponding part
of equations.
But let us stress the fact that this depends on the coordinate
prescription given.

Such coordinates exist for $g=0$ (the quasi-global coordinate $z$)
and for $g=1$ (the coordinate on the simply-connected covering).
In Section \kzrat\ and Section \kzell\ we are exactly dealing with this
situation.
For higher genus all uniformisations have certain disadvantages.
We could either realize the curve as upper half-plane modulo
a Fuchsian group or embed it into its
Jacobian torus.
In the first case we do not have a nice behaviour under deformation
of the complex structure. In the  second case we have multi-dimensional
coordinates.

\vskip 1cm
\head
  \kzrat . The case of lower genus: $g=0$
\endhead
\vskip 0.2cm
\numsec=\kzrat   
\numfor=1          

 Let us show how to obtain the original Knizhnik-Zamolodchikov
 equations \equ(KZeqi) from \equ(kze) for $g=0$.
 Let  $z_i$ ($i=1,\ldots ,N$) be the $N$ moving points
 and fix the reference point $z_\infty$
 to be $\infty$.
 For definiteness take as representation sheaf a Verma module sheaf
 as introduced in \equ(vermash), \equ(vermam).
 We set   $ \     \ai:=\prod\limits_{l=1\atop l\ne i}^N (z_i-z_l)^{-1}\ $
 for $i=1,\ldots,N$.
 For every point   $z_i$ ($i=1,\ldots ,N$) the
 KN basis elements of degree $(-1)$ for the quadratic differentials,
 resp. for the  vector  fields  are given as
$$\Omega^i(z):= \Omega^{-1,i}=\frac{\text{d}z^2}{z-z_i},\quad
\text{resp\.}\quad
e_i(z):=e_{-1,i}= \big(\ai\prod_{j\ne i}(z-z_j)\big)
\frac{\partial}{\partial
 z}\ .\Eq(bqvf)$$
 The vector
 field $e_i$  evaluates  to $\frac{\partial}{\partial z}$  at
 the point $z_i$ and vanishes at all other points $z_j,$
 $j\ne i.$ Therefore $e_i$ corresponds to the  basic
 direction  $\partial_i$ on the configuration space
 which is responsible for moving
 the point $z_i$. On the other
 hand, $e_i$ is exactly the \KN\ dual vector field to the
 quadratic differential $\Omega^i$. This  follows from
 calculating the residues, see \equ(edu) .

 The coefficients $l^{(m,i)(n,j)}_k$ are given by
   $$ l^{(m,i)(n,j)}_k = \cintt\omega^{(m,i)}\omega^{(n,j)}e_k\ .
                                                \Eq(coef)
   $$
Recall from \equ(ebgo)
   $$ \omega^{(m,i)} = \frac{\ai^{-m}\;dz}{(z-z_i)\prod\limits_{s=1}^N
                       (z-z_s)^m}\ .
   $$
 Therefore the integrand in \equ(coef) equals to
 $$ \frac{\ai^{-m}\aj^{-n}\ak\; dz}{(z-z_i)(z-z_j)(z-z_k)\prod\limits_{s=1}^N
    (z-z_s)^{m+n-1}} \ .
                                              \Eq(plotn)
 $$
 The coefficients \equ(coef) can be obtained by summation of the
 residues at the points $z_1,\ldots,z_N$ (or  alternatively
 by the negative of the residue at the
 point $z_\infty$). Note that if ($m+n\le -2$) or ($m+n=-1$ and not
 $i=j=k$) then
 all residues at the points $z_1,\ldots,z_N$ vanish. If $m+n>0$
 the residue vanishes at $z_\infty$.
 Here we consider only the case  when
 $\Phi$ in \equ(kze) is a vector which is annihilated by the
 subalgebra $\gh_+$; note that the KZ equations were originally
 obtained under the same assumption in \cite\KnZ.  But, if $m+n=0$
 and $m,n\ne0$ then either $m$ or $n$ is positive and hence
 $\nord{u(m,i)u(n,j)}\Phi =0$, because by the normal ordering the
 elements of positive degree will appear on the right.  The  nonzero
 coefficients in \equ(kze) for a given $k=1,\ldots,N$ are as follows:
$$
\gathered
l^{(0,i)(0,i)}_k=\frac {\ai^{-1}\ak}{z_i-z_k},\  (i\ne k),\quad
 l^{(0,i)(0,k)}_k=l^{(0,k)(0,i)}_k =\frac
 {1}{z_k-z_i},\  (i\ne k),
\\
 l^{(0,k)(0,k)}_k
 =\sum\limits_{i\ne k}\frac {1}{z_k-z_i},\qquad
 l^{(-1,k)(0,k)}_k=l^{(0,k)(-1,k)}_k= \ak^2\ .
\endgathered
\Eq(lequ)
$$
In the remainder of this section we will show that after
modifying the vector fields $e_k$ by adding vertical vector fields
(i.e. vector fields which have zeros at all the points
$z_i$ and hence are not moving these points),
applying a certain factorization process, and calculating the structure
constants with respect to this basis,
that all
coefficients
can be
eliminated except the one with $m=n=0$ and
nonequal
upper indices.
 Hence the equation \equ(kze)  (now with respect to the
  modified basis)  will have the following form:
 $$\bigg(\partial_i-\frac {1}{\ce+\ka}
     \sum\limits_{j\ne i}
     \frac{\nord{u(0,i)u(0,j)}+\nord{u(0,j)u(0,i)}}{z_i-z_j}
     \bigg)\Phi =0\ ,\quad i=1,\ldots ,N.      \Eq(KZres)
  $$
 Note that the coefficients $l^{(0,i)(0,j)}_k$ ($i\ne
 j$) can also be obtained from the expansion
     $$\split
     \omega^{0,i}(z)\omega^{0,j}(z) &=
     \frac{\text{d}z^2}{(z-z_i)(z-z_j)}=
     \frac{1}{z_i-z_j}\bigg(
      \frac{\text{d}z^2}{z-z_i}-
      \frac{\text{d}z^2}{z-z_j}
     \bigg) \\ &=
     \frac{1}{z_i-z_j}(\Omega^i-\Omega^j)\,\,(i\ne j)
     \endsplit ,
   $$
which implies
$$l^{(0,i)(0,j)}_k=\cases  0, &k\ne i\,\,\text{and}\,\,k\ne j\\
            \frac{1}{z_k-z_j}, &k=i\\
            -\frac{1}{z_i-z_k}, &k=j\ .
           \endcases
$$
 Let us improve the vector fields $e_k, k=1,\ldots, N$ in order to
 eliminate the contribution of the terms with $i=j$, $m=n=0$
 (without changing the expressions for $i\ne j$).
  We set ${\tOmega}^i :=
 \omega^{0,i}\omega^{0,i}=\frac{dz^2}{(z-z_i)^2}$. Let us denote the
 \KN\ pairing \equ(knpair) between 2-differentials and vector fields by
 the angle   brackets $\langle\cdot ,\cdot\rangle$. Then
  $$ \langle{\tOmega}^i ,e_k\rangle = \frac{1}{z_k-z_i}
     \prod\limits_{s\ne i,k} \frac{z_i-z_s}{z_k-z_s},
     \,\,(k\ne i)\ ,\qquad
     \langle{\tOmega}^k ,e_k\rangle = \sum\limits_{s\ne k}
     \frac{1}{z_k-z_s}.
                                                  \Eq(pair)
  $$
Let us pass from the set of vector fields $\{e_k\, |\, k=1,\ldots
,N\}$ to the set of vector fields $\{e_k^\prime\, |\, k=1,\ldots
,N\}$ where
  $$ e_k^\prime := e_k + \sum\limits_{i=1}^N \l_{ki}E_i,
  \qquad \l_{ki}\in\C\ ,
                                                      \Eq(eprime)
  $$
and $\{E_s\, |\, s=1,\ldots ,N\}$
are the \KN\ basis elements  of degree zero for the vector fields, i.e
  $$ E_{i} :=e_{0,i}=
  (z-z_i)\prod\limits_{s\ne i} \frac{(z-z_s)^2}{(z_i-z_s)^2}
                   \frac{\partial}{\partial z}.       \Eq(corr)
  $$
\proclaim{Lemma \kzrat .1} (a) $E_i$ is a vertical vector field
 ($i=1,\ldots ,N$).
 \nl
 (b) $\langle \tOmega^i, E_j\rangle = \d_{ij}$ ($i,j=1,\ldots ,N$)
 \nl
 (c) $\langle \Omega^i, E_j\rangle = 0$ ($i,j=1,\ldots ,N$)
\endproclaim
\demo{Proof} (a) follows from the fact that the vector fields $E_k$
have zeros at all the points $\{z_s\, |\, s=1,\ldots ,N\}$.

(b) By definition
   $$ \langle \tOmega^i, E_j\rangle = \sum\limits_{p=1}^N res_{z_p}
      \frac{(z-z_j)}{(z-z_i)^2} \prod\limits_{s\ne j}
      \frac{(z-z_s)^2}{(z_j-z_s)^2} \, dz.
   $$
If $i\ne j$ then the 1-form on the right hand side of the latter
relation is holomorphic and all residues equal to zero. If $i=j$ then
all the residues are zero except at the point $z_i$ and there
it is equal to 1.

(c) All 1-forms $\Omega^iE_j$, $i,j=1,\ldots,N$ are holomorphic which
proves the claim, see also \equ(edu).
\qed
\enddemo
By Lemma \kzrat .1(a)  the vector fields $\{e_k^\prime\, |\,
k=1,\ldots,N\}$ correspond
to the basic infinitesimal deformations
$\partial_i$ too,
as well as the vector fields $\{e_k\, |\, k=1,\ldots,N\}$
do.  By Lemma \kzrat .1(c) the coefficients $l^{(m,i)(n,j)}_k$ ($i\ne
j$) remain the same under replacing the set of vector fields $\{e_k\,
|\, k=1,\ldots,N\}$ with the set $\{e_k^\prime\, |\, k=1,\ldots,N\}$.

Now let us find the coefficients $\l_{ki}$ in \equ(eprime) in such a
way that
  $$ \langle e_k^\prime ,\tOmega^i\rangle = 0,
     \,\,\,\, i,k=1,\ldots,N.                  \Eq(noCont)
  $$
This  means that $l^{(0,i)(0,i)}_k=0$, $(i=1,\ldots,N)$ if
these coefficients are calculated
with respect to the  vector fields $\{e_k^\prime\, |\,
k=1,\ldots,N\}$.  Recall that $\tOmega^i=\omega^{0,i}\omega^{0,i}$.

The equation \equ(noCont) means
  $$ \langle e_k +\sum\limits_{s=1}^N \l_{ks}E_s,\tOmega^i\rangle =0\ ,
     \qquad i,k=1,\ldots,N.                  \Eq(noCont1)
  $$
By Lemma \kzrat .1(b) this is equivalent to $\langle e_k,
\tOmega^i\rangle + \l_{ki} = 0$.
Hence \equ(noCont) if and only if
$\  \l_{ki}=-\langle e_k,\tOmega^i\rangle$.

Taking into account \equ(pair) this
enables us to give explicit expressions for the adjusted $e_k^\prime$'s:
  $$ e_k^\prime = e_k + \sum\limits_{i\ne k}
     \frac{\prod\limits_{s=1}^N (z-z_s)^2}{z_k-z_i} \cdot\ak
     \left( \frac{\ai}{z-z_i} - \frac{\ak}{z-z_k}\right)
     \frac{\partial}{\partial z},
  $$
where $i,k=1,\ldots,N$.
Note that there might occur now
nonvanishing coefficients  with $m+n=1$. But again under the
assumption that $\gh_+$ annihilates $\Phi$ and with the normal ordering
the corresponding operator terms will not contribute to the final equation.

At last the term
$l^{(0,k)(-1,k)}_k(\nord{u(0,k)u(-1,k)}+\nord{u(-1,k)u(0,k)})$
of the Knizhnik-Zamolodchikov equation can be eliminated by
either considering elements  only of degree zero
in the result or by
passing to the  quotient sheaf
$\widehat{V}_{\l,\gamma,\d}/\gb_{-}^*\widehat{V}_{\l,\gamma,\d}$
(see \equ(cbdef)) .
Of course, one has to assume
$\nord{u(0,k)u(-1,k)}\,=\,\nord{u(-1,k)u(0,k)}\,=u(-1,k)u(0,k)$ for
doing that.
But this is true in the standard normal ordering \equ(normst).

>From \equ(KZres) the exact form of \equ(KZeqi) can be obtained
by taking $\Phi$ to be from the the $\g$-module
$V_{\lambda,\gamma}$ (see Proposition \kzrep.4) which is
a tensor product of the individual representations $V_{\lambda_i}$.
If we assume that the $u_a$ are a selfdual basis of $\g$ and take the
standard normal ordering then
$$
u_a(0,i)\Phi=t^a_i\Phi,\quad
\text{and further}
\quad
u_a(0,j)u_a(0,i)\Phi=t_j^at^a_i\,\Phi=t_i^at^a_j\,\Phi\quad \text{for }
 j\ne i\ .
$$
Hence,
$$
\bigg(\frac{\partial}{\partial z_i}-
     \frac {2}{\ce+\ka}\sum\limits_{j\ne i}\frac{t^a_it^a_j}{z_i-z_j}\bigg)
    \Phi=0\ ,
                                    \qquad i=1,\ldots,N\ .
$$



\vskip 1cm
\head
  \kzell . The case of lower genus: $g=1$
\endhead
\vskip 0.2cm
\numsec=\kzell   
\numfor=1          

The purpose of this section is to obtain explicit expressions for
the coefficients of the KZ equations via the Weierstrass $\s$-function.

{\bf 1.} Let us take the following set of vector fields corresponding to
the motion of the points $z_1,\ldots,z_N$: $e_k(z)=A_{0,k}(z)
\frac{\partial}{\partial z}$ ($k=1,\ldots,N$, see Section~\kzapp), or
explicitly
  $$ e_k(z)= \prod\limits_{s\ne k} \frac{\s(z-z_s)}{\s(z_k-z_s)}
     \cdot
     \frac{\s(z_k-z_0)^N}{\s(z-z_0)^N} \cdot
      \frac{\s(z+\sum_{s\ne k}z_s-Nz_0)}
           {\s(\sum_{s=1}^N z_s-Nz_0)}
      \frac{\partial}{\partial z}.               \Eq(fieldk)
  $$
Here $z_0$ is the fixed reference point.
Our first goal is to find the contribution of certain terms of the
form $l^{(m,i)(n,j)}_k\nord{u(m,i)u(n,j)}$ into the KZ equations.  Note
that by the duality relations $\w^{m,i} = A_{-m-1,i}dz$.  We will consider
several cases.
As in genus zero let
$\Phi$  be an element of the representation
space
which is
annihilated by the subalgebra
$\gh_+$.

{\bf 1.1.} $m\ne 0$, $n\ne 0$.
\nl
Define $\tw^{m,i}$, $\tw^{n,j}$, ${\tilde e}_k$ by the relations
$\w^{m,i}=\s(z-z_i)^{-1}\tw^{m,i}$,
$\w^{n,j}=\s(z-z_j)^{-1}\tw^{n,j}$,
$e_k=\s(z-z_k)^{-1}{\tilde e}_k$.
Then $ord_{z_s}\tw^{m,i}=-m$,
 $ord_{z_s}\tw^{n,j}=-n$,  $ord_{z_s}{\tilde e}_k=1$
($s=1,\ldots,N$), $ord_{z_0}\tw^{m,i}=Nm$, $ord_{z_0}\tw^{n,j}=Nn$,
$ord_{z_0}{\tilde e}_k=-N$.
Set $\tw^{(m,i)(n,j)}_k:=\tw^{m,i}\tw^{n,j}{\tilde e}_k$. Then
one has
  $$\w^{m,i}\w^{n,j}e_k =
    \frac{\tw^{(m,i)(n,j)}_k}{\s(z-z_i)\s(z-z_j)\s(z-z_k)}.
                                              \Eq(c1.1)
  $$
and $ord_{z_s}\tw^{(m,i)(n,j)}_k = -m-n+1$ for all $s=1,\ldots,N$.
If $m+n\ge 1$ then \equ(c1.1) is holomorphic at $z_0$. If $m+n=0$
then either $m>0$ or $n>0$. In both cases
 $\nord{u(m,i)u(n,j)}\Phi =0$.
By the same reason for the case
$m+n=-1$ one would have $m=0$, $n=-1$ (or vice versa).
But in this subsection we assume
$m\ne 0$, $n\ne 0$, so this case will not appear.
It remains $m+n\le -2$. Then for the order of the numerator of
 \equ(c1.1) one  has $ord_{z_s}\tw^{(m,i)(n,j)}_k=-m-n+1\ge 3$
 ($s=1,\ldots,N$), hence the 1-form is holomorphic at all the marked
 points even if some or all of the points $z_i$, $z_j$, $z_k$
coincide.

Hence the case under
consideration does not contribute to the KZ equations.

{\bf 1.2.} $m=0$, $n\ne 0$ (or vice versa).
\nl
By (\kzapp.4) $\w^{0,i}=\hw^{0,i}-\sum\limits_{s=1}^N\gamma_{i,s}\w^{-1,s}$
 where $\hw^{0,i}=A^\prime_{-1,i}dz$ (see (\kzapp.7)). The order of the sum
 on the right hand of the latter relation at any moving point is determined
 by the order of its first summand because all other terms have some bigger
 orders. If $n>0$ then $\nord{u(0,i)u(n,j)}\Phi =0$.
 Using the definition of $e_k$ and
 $\w^{m,i}=A_{-m-1,i}dz$ one obtains $\hw^{0,i}\w^{n,j}e_k=
 A^\prime_{-1,i}A_{-1-n,j}A_{0,k}\, dz$. By the explicit
 expressions \equ(fieldk) , (\kzapp.6)-(\kzapp.8) one obtains the following:
 if $n\le -2$ then $-n+1\ge 3$ and  $\hw^{0,i}\w^{n,j}e_k$
 is holomorphic at any
 marked point. As $n\ne 0$ it remains only to consider the case
 $n=-1$. In this case  $\hw^{0,i}\w^{n,j}e_k$ has a pole
 of order at most one at the moving points. And this will occur
 only if $i=j=k$.
 In particular, the terms
 with $\gamma^\prime$s are holomorphic and do not contribute to the
 result. An easy calculation of  the
 residues gives
  $$ l^{(0,k)(-1,k)}_k= 1\,\,\, (k=1,\ldots, N).
                                                     \teq(coef0)
  $$
  In the same way we obtain the symmetric expression
 $l^{(-1,k)(0,k)}_k=1$.

{\bf 1.3.} $m=n=0$.
\nl
To find $l_k^{(0,i)(0,j)}$ one has to consider 1-forms $\w^{0,i}
\w^{0,j}e_k$ where $\w^{0,i}=\hw^{0,i}-\sum\limits_{s=1}^N \gamma_{is}
\w^{-1,s}$, $\w^{0,j}=\hw^{0,j}-\sum\limits_{r=1}^N \gamma_{jr}\w^{-1,r}$.
 The terms $\w^{-1,s}\w^{-1,r}e_k$ are holomorphic
because all their factors are. A term of the form $\w^{-1,s}\hw^{0,j}e_k$
 can have a pole (which is necessarily of order 1)
 only if $s=j=k$ (and is holomorphic otherwise).
So one has
$$l_k^{(0,i)(0,j)}=\langle\hw^{0,i}\hw^{0,j},e_k\rangle -
 (\gamma_{jk}\delta_i^k+\gamma_{ik}\delta_j^k)
\langle \w^{-1,k}\hw^{0,k},e_k \rangle\ .$$
 The second scalar product
 was  found already and equals $l^{(0,k)(-1,k)}_k$ which is
 simply $1$.  For the first term one has
 $\hw^{0,i}\hw^{0,j}e_k=A^\prime_{-1,i}A^\prime_{-1,j}A_{0,k}\, dz$.
In case $i\ne j$ this 1-form has a residue at the point $z_i$ if
$i=k$ or at the point $z_j$ if $j=k$ (if $i,j,k$ are pairwise
different then the 1-form is holomorphic). For $i=k\ne j$ one has
  $$ l^{(0,k)(0,j)}_k  = \frac{1}{\s(z_k-z_j)}
     \frac{\s(z_k-w_1)\s(z_k-w_2)\s(z_j-z_0)}
         {\s(z_j-w_1)\s(z_j-w_2)\s(z_k-z_0)} - \gamma_{jk}
                                                \Eq(coef1)
  $$
and an analogous expression for $i\ne k=j$.
There is $w_1+w_2 = z_i+z_0$  in \equ(coef1) and $w_1$ is the same as
in (\kzapp.8).
Note that
in  case $i=k\ne j$
this coefficient comes with the operator-valued factor
$\nord{u(0,k)u(0,j)}$   and with the factor
$\nord{u(0,i)u(0,k)}$ in case $i\ne k=j$.

In the same way as in case
$g=0$ we can annihilate the contribution of $\w^{0,i}\w^{0,i}$ by
adding certain vector fields of degree $0$ to $e_k$.
So we come to the following
form of the KZ equations corresponding to  moving the points ($k=1,\ldots,N$):
  $$ \split
     \partial_k\Phi  & -\frac {1}{\ce+\ka}
      \sum\limits_{i\ne k} l^{(0,k)(0,i)}_k
     \left( \nord{u(0,k)u(0,i)} +\nord{u(0,i)u(0,k)} \right)
     \Phi       \\
     &  -\frac {1}{\ce+\ka}   \sum\limits_{i=1}^N
     \left( \nord{u(0,i)u(-1,i)} +\nord{u(-1,i)u(0,i)} \right)
     \Phi       =0\ .
     \endsplit
                                                      \teq(mvpnt)
  $$

\noindent
{\bf 2.} Now let us consider the KZ equation corresponding to the
deformation of the complex structure (i.e. to a change of the
moduli parameter). Consider the following vector field:
  $$ e_0= \s(z-E)^{N+1} \s(z-z_0)^{-1}
     \prod\limits_{s=1}^N \s(z-z_s)^{-1}
     \frac{\partial}{\partial z}\ ,
  $$
where $E=(N+1)^{-1}(z_0+z_1+... +z_N)$.

The vector field $e_0$ has simple poles at all the points
$z_1,\ldots,z_N,z_0$.
It follows from the following lemma that
the corresponding tangent
vector
on the
moduli space is non-trivial and that it is independent
of the ones generated by the vector fields $e_1,\ldots,e_N$.
Hence, it corresponds indeed to a deformation of the complex structure.
In the lemma vertical vector fields are considered
with respect to the moduli space $\MgN$.
\proclaim{Lemma \kzell .2} For  points $z_1,\ldots,z_n,z_0$ in
generic position the vector field $e_0$ cannot be expressed as a
linear combination of the vector fields $e_1,\ldots,e_N$ and vertical
vector fields.
\endproclaim
\demo{Proof} As the vector fields $e_1,\ldots,e_N$ are regular at
the points $z_1,\ldots,z_N$ their residues at the point $z_0$ are
zero (on the elliptic curve one can speak about residues of a
vector field). All vertical vector fields have zero residues at
$z_0$ too. This is evident for the ones which are regular at the
point $z_0$. As for the ones which have zeroes at the points
$z_1,\ldots,z_N$ they also have zero residues at the point $z_0$
because the sum of all their residues vanishes. But the vector field
$e_0$ has a simple pole at $z_0$ with nonzero residue for a generic
point of the moduli space under consideration. This proves the lemma.
\qed
\enddemo
Let us consider several cases as earlier.

{\bf 2.1.} $m,n\ne 0$.
\nl
Define $\tw_0^{(m,i)(n,j)}$ from the relation
  $$ \w^{m,i}\w^{n,j}e_0=
     \frac{\tw_0^{(m,i)(n,j)}}{\s(z-z_i)\s(z-z_j)},
                                               \Eq(c2.1)
  $$
Then $ord_{z_s}{\Tilde\w}_0^{(m,i)(n,j)} = -m-n-1$ $(s=1,\ldots,N)$.
As the order of the denominator of \equ(c2.1) could be
 at most two the
 1-form \equ(c2.1) is holomorphic at the
points $z_1,\ldots,z_N$ as soon as $-m-n-1\ge 2$,
i.e. if $m+n\le -3$. If $m>0$ or $n>0$ then $\nord{u(m,i)u(n,j)}\Phi =0$.
So one has either $m+n=-1$ or $m+n=-2$. As
$m,n\ne 0$ the first case does not occur and in the second case
only  $m=n=-1$ remains.
 Also $i=j$ in this case because otherwise \equ(c2.1) is
holomorphic. Using explicit formulas one obtains
  $$ l^{(-1,i)(-1,i)}_0= \s(z_i-E)^{N+1} \s(z_i-z_0)^{-1}
     \prod\limits_{s\ne i}\s(z_i-z_s)^{-1}
                                                  \Eq(new1)
  $$
for $i=1,\ldots,N$.

{\bf 2.2.} $m=0, n\ne 0$ (or vice versa).
\nl
In this case one has to put $m=0$
 in \equ(c2.1). Then $ord_{z_s}\tw^{(0,i)(n,j)}=-(n+1)$ ($s=1,\ldots,N$).
 If $n>0$ then $\nord{u(0,i)u(n,j)}\Phi =0$. If $-(n+1)\ge 2$, i.e.
 $n\le -3$ the 1-form \equ(c2.1) is holomorphic. Hence one has only to
 consider $n=-1$ and $n=-2$.

If $n=-1$ then $\w^{0,i}\w^{-1,j}e_0=\hw^{0,i}\w^{-1,j}e_0-
\sum\limits_{s=1}^N \gamma_{i,s}\w^{-1,s}\w^{-1,j}e_0$.
If $s\ne j$ then the 1-forms $\w^{-1,s}\w^{-1,j}e_0$ are
 holomorphic at all  moving points (see above).
If $s=j$ then the contribution
 of this 1-form equals to $l^{(-1,j)(-1,j)}_0$ (see Subsection 2.1).
 Furthermore
  $$ \hw^{0,i}\w^{-1,j}e_0= \frac{F(z)}{\s(z-z_i)\s(z-z_j)}dz,
  $$
where
  $$ F(z):= A^\prime_{-1,i}(z)A_{0,j}(z) \s(z-z_0)^{-1} \s(z-E)^{N+1}
     \s(z-z_i)\prod\limits_{s\ne j}\s(z-z_s)^{-1}\ ,
  $$
and $A^\prime_{-1,i}(z)$, $A_{0,j}(z)$ are introduced in
Section~\kzapp .  As follows from
Section ~\kzapp
that $F(z)$ has order zero
at the points
$z_1,\ldots,z_N$.
If $i\ne j$ then there are two residues at the point $z_i$ and at the
point $z_j$. So we obtain
  $$ l_0^{(0,i)(-1,j)} =
     \frac{F(z_i)-F(z_j)}{\s(z_i-z_j)}
      -\gamma_{ij}l^{(-1,j)(-1,j)}_0 \qquad\text{for}\quad i\ne j .
                                            \teq(co2.2)
     $$
This residue contributes to the KZ equation with the operator-valued
coefficient
\nl
$\nord{u(0,i)u(-1,j)}$. Furthermore
$l_0^{(-1,i)(0,j)}=l_0^{(0,j)(-1,i)}$ but it comes with the
operator-valued coefficient $\nord{u(-1,i)u(0,j)}$.
In the case $i=j$  poles of 2-d order arise and one has
  $$l_0^{(0,i)(-1,i)} =
    \left.
      \frac{d}{dz}\left( (z-z_i)^2 F(z)\right)
    \right|_{z=z_i}-
      \gamma_{ii}l^{(-1,i)(-1,i)}_0.
  $$

Consider now   $n=-2$.
 Then the order of the numerator of \equ(c2.1) equals 1 at all the
 moving points. So the residue can be nonzero only if $i=j$. The
additional terms $\gamma_{is}\w^{-1,s}\w^{-2,i}e_0$ are holomorphic
 and will not contribute into the result. By explicit formulas one has
   $$ l^{(0,i)(-2,i)}_0=
          \s(z_i-E)^{N+1} \s(z_i-z_0)^{-1}
          \prod\limits_{s\ne i}\s(z_i-z_s)^{-1}\ .
                                                \teq(new2)
    $$

{\bf 2.3.} $m=n=0$.
\nl
One has
 $$\w^{0,i}\w^{0,j}e_0=\hw^{0,i}\hw^{0,j}e_0-\sum_{s=1}^N \gamma_{is}
 \w^{-1,s}\hw^{0,j}e_0-\sum_{r=1}^N \gamma_{jr} \w^{-1,r}\hw^{0,i}e_0+
 \sum_{s,r=1}^N \gamma_{is}\gamma_{jr} \w^{-1,s}\w^{-1,r}e_0.
 $$
Contributions of the
 terms with $\gamma^\prime$s can be expressed via the coefficients calculated
in the Subsections 2.1 and 2.2. For the first term of the above sum one has
$\hw^{0,i}\hw^{0,j}e_0 =F_{ij}(z)\, dz$, where
  $$\align
    F_{ij}(z)= C_{-1,i}^{-1}C_{-1,j}^{-1}
    \frac{\s(z-w_1)^2\s(z-w_2)^2}{\s(z-z_i)\s(z-z_j)}
          & \prod\limits_{s=1}^N \s(z-z_s)^{-1} \teq(c2.3) \\
   \times & \s(z-E)^{N+1}\s(z-z_0)^{-3}
    \endalign
  $$
and the coefficients are given by (\kzapp.8).
It has poles of  2-d order and if $i=j$ then even of
3-d order at the points $z_i, z_j$.
The corresponding coefficients look as follows:
  $$l^{(0,i)(0,j)}_0=\langle \hw^{0,i}\hw^{0,j},e_0\rangle  -\sum
    \limits_{s=1}^N (\gamma_{is}l^{(0,j)(-1,s)}_0+\gamma_{js}
    l^{(0,i)(-1,s)}_0+\gamma_{is}\gamma_{js}l^{(-1,s)(-1,s)}_0)
                                                       \Eq(co2.3)
  $$
where
$
\langle \hw^{0,i}\hw^{0,j},e_0\rangle =
\left.\frac{d}{dz}\left( (z-z_i)^2 F_{ij}(z)\right)\right|_{z=z_i}+
\left.\frac{d}{dz}\left( (z-z_j)^2 F_{ij}(z)\right)\right|_{z=z_j}
$ ($i\ne j$),
$\langle \hw^{0,i}\hw^{0,i},e_0\rangle =  \frac{1}{2}
\left.\frac{d^2}{dz^2}\left( (z-z_i)^3 F_{ij}(z)\right)\right|_{z=z_i} $.
All the terms of the KZ equation which contain $u(n,j)$, $n<-1$ can be
eliminated by factorization over $\gb_{-}^*$ as in Section~\kzrat\ and we
obtain the equation corresponding to the deformation of
the complex structure  in the
following form:
  $$\align
    \partial_0\Phi
   & -\frac{1}{\ce+\ka} \sum\limits_{i,j=1}^N
    l^{(0,i)(0,j)}_0 \nord{u(0,i)u(0,j)}\Phi    \\
    & -\frac{1}{\ce+\ka}
    \sum\limits_{i,j=1}^N l_0^{(0,i)(-1,j)}
    (\nord{u(0,i)u(-1,j)}+  \nord{u(-1,j)u(0,i)})\Phi   \teq(KZ0)  \\
    & -\frac{1}{\ce+\ka}  \sum\limits_{i=1}^N
    l^{(-1,i)(-1,i)}_0 \nord{u(-1,i)u(-1,i)}\Phi
     \quad=\quad 0\ .
    \endalign
  $$


\vskip 1cm
\head
  \kzapp.  Appendix:  KN Basis for the elliptic case
\endhead
\vskip 0.2cm
\numsec=\kzapp   

\def\kznum{\kzapp}
\numfor=1      
%
%
In the elliptic case (i.e. $g=1$) the canonical bundle $\K$
and hence all its powers are trivial.
This implies that for all weights $\l$ there is a change of
the prescription for certain basis elements necessary.
The adopted prescription is given in this appendix.
Furthermore, explicit expressions in terms of
Weierstra\ss-$\sigma$ function are given.

By the triviality of the
canonical bundle
we have the following relation between the basis elements of $\Fl$
$$
f_{n,p}^\l=A_{n-\l,p}dz^\l\ .\Eq(elre)
$$
In particular,  by fixing  a basis of $\Fn0$ we obtain
one  for every $\Fl$.
\nl
The standard prescription for the $(N,1)$ situation
with incoming points $\{P_1,P_2,\ldots,P_N\}$ and outgoing point
$\{\Pif\}$ is
$$
\ord_{P_p}(A_{n,p})=n,\quad
\ord_{P_i}(A_{n,p})=n+1,\ \text{for\ }i\ne p,\qquad
\ord_{\Pif}(A_{n,p})=-N(n+1)\ .
$$
For the $(1,1)$ situation we set as above $\ A_0=1\ $ and
$A_{-1}'dz=\rho=(\omega^0)'$.
If we want to have duality we have to modify these elements.
With
$$
c=\cintt {A_{-1}'}^2dz\Eq(elcaa)
$$
we set
$$
A_{-1}:=A_{-1}'-\frac c2 A_0 \Eq(elca)
$$
and obtain the required
relation
e.g. $\langle A_{-1},\omega^0\rangle=0$.
Note that $\omega^0=A_{-1}dz$
is not obliged to have purely imaginary periods.

For the $(N,1)$ situation with $N>1$ we have the standard
prescription for $A_{0,p}$ and hence for $\omega^{-1,p}$
and the following modified prescription for
$n=-1$
$$
\ord_{P_p}(A_{-1,p}')=-1,\quad
\ord_{P_i}(A_{-1,p}')=0,\ \text{for\ } i\ne p, \quad
\ord_{\Pif}(A_{-1,p}')=-1\ .
$$
By adding certain linear combinations of $A_{0,s}$ to
$A_{-1,p}'$ we obtain again a basis which fulfils the
duality.
In more detail, set
$$
\gamma_{r,p}:=\frac 12\cintt A_{-1,r}'A_{-1,p}'dz\ . \Eq(elcb)
$$
The adjusted basis elements are given as
$$
A_{-1,r}:=A_{-1,r}'-\sum_{s=1}^N\gamma_{r,s}A_{0,s}\ .\Eq(elcc)
$$

Explicit representations of the basis elements for $g=1$
in terms of Weierstra\ss\
$\sigma$-function are as
follows \cite{\rSLb}:

Let the torus be given as
$\ T=\C\bmod L\ $ with the normalized lattice $L=\langle 1,\tau\rangle{}_\Z$.
Let $z_i\in\C$, $i=1,..,N$
 be fixed lifts of the points $P_i\in T$  to $\C$, i.e.
$z_i\bmod L=P_i$ and $z_0$ be a fixed lift of the point $\Pif\in T$.
Let $n\in \Z$ be  fixed.
For $p=1,\ldots, N$ set
$$b_{n,p}:=b_{n,p}(z_1,\ldots,z_N,z_0)
:=-(n+1)\sum_{i=1}^Nz_i+z_p+N(n+1)z_0\ .
$$
If $n\ne -1$ for $N>1$ or $n\ne 0,-1$ for $N=1$
then
for  generic choices of the points $z_i$ and the point $z_0$ the
$b_{n,p}$
does
not coincide $\bmod L$ with them, even if
we move them locally.
But note that the points $b_{n,p}$ depend on the variation of the points $z_i$
and $z_0$.
At the point $P_i$ we choose as coordinate $z_i=(z-z_i), i=0,1,\ldots, N$.
\proclaim{Proposition \kznum.1}\cite{\rSDiss, p.53}
Let the  points $P_i, i=1,\ldots, N$ and $\Pif$ be generic.
Then the basis element
$A_{n,p}$ for ($n\ne 0,-1$ and $N=1$) and for ($n\ne -1$
and $N>1$) is given as
$$
A_{n,p}(z)=C^{-1}_{n,p}\bigg(\prod_{i=1}^N\sigma(z-z_i)\bigg)^{n+1}
\sigma(z-z_p)^{-1}\sigma(z-z_0)^{-N(n+1)}\sigma(z-b_{n,p})\ ,\Eq(ebeg)
$$
with
$$
C_{n,p}:=\bigg(\prod_{i\ne p}\sigma(z_p-z_i)\bigg)^{n+1}
\sigma(z_p-z_0)^{-N(n+1)}\sigma(z_p-b_{n,p})
\ .
$$
In particular, it is a well-defined function on the torus.
\endproclaim
The point $b_{n,p}\bmod L$ is the additional zero of the
basis element $A_{n,p}$.
Note that the constants $C_{n,p}$ depend also on the points
$z_0,\ldots,z_N$.

For the exceptional values one obtains the following
results.
\nl
For $n=0$: If $N=1$ one sets  $A_0:=1$.
For $N>1$ the above formula is already the correct formula.
It specialises to
$$
A_{0,p}=C^{-1}_{0,p}\bigg(\prod_{i\ne p}\sigma(z-z_i)\bigg)
\sigma(z-z_0)^{-N}
\sigma(z+\sum_{i\ne p}z_i-Nz_0)\Eq(ebeo)
$$
with
$$
C_{0,p}:=\bigg(\prod_{i\ne p}\sigma(z_p-z_i)\bigg)\sigma(z_p-z_0)^{-N}
\sigma(z_p+\sum_{i\ne p}z_i-Nz_0)\ .
$$

For $n=-1$:
In this case the order at the point $\Pif$ is set to $-1$.
We have to find elements $w_1,w_2\in\C$ such that $w_1+w_2=z_p+z_0$
but $w_1,w_2\ne z_i\bmod L,\ \i=0,\ldots, N$.
It is possible to choose $w_1$ generic and this will fix $w_2$.
\proclaim{Proposition \kznum.2} \cite{\rSDiss}
For generic choice of points we have
$$
A_{-1,p}=C^{-1}_{-1,p}\sigma(z-z_p)^{-1}\sigma(z-z_0)^{-1}\sigma(z-w_1)
\sigma(z-(z_p+z_0-\w_1))\Eq(ebem)
$$
with
$$
C_{-1,p}:=\sigma(z_p-z_0)^{-1}\sigma(z_p-w_1)
\sigma(\w_1-z_0)\ .
$$
\endproclaim
The arbitrariness of the choice in $w_1$ comes from the
fact that the element is only fixed
by the orders at $P_i$ and $\Pif$ up to
addition of an arbitrary constant.

Now $(\omega^{0,p})'=A_{-1,p}'dz$ and $\omega^{-1,p}=A_{0,p}dz$.
But note that the duality condition is not yet fulfilled.
In particular we do not have
$\langle \omega^{0,p},A_{-1,s}'\rangle=0$.
To obtain this condition we have to add  combinations of the
elements $A_{0,r}$ to the $A_{-1,s}'$, as it is expressed by
\equ(elcb), \equ(elcc).

%
%
\define\PL{Phys\. Lett\. B}
\define\NP{Nucl\. Phys\. B}
\define\LMP{Lett\. Math\. Phys\. }

\define\CMPP{Commun\. Math\. Phys\. }
\define\JMP{J\.  Math\. Phys\. }

\define\FA{Funktional Anal\. i\. Prilozhen\.}

{
}

\enddocument
\bye